\input amstex
\documentstyle{amsppt}
\magnification=1200
\NoBlackBoxes
\TagsOnRight
\NoRunningHeads
\vsize=22 truecm
\hsize=16 truecm

\def \real{{\Bbb R}}
\def \complex{{\Bbb C}}

\def \integer{{\Bbb Z}}

\redefine \natural{{\Bbb N}}
\def \sp {{ \text{sp}\, }}
\def \var {{ \text{var}\, }}
\def \ess {{ \text{ess}\, }}
\def \Id {{\text{Id}\, }}
\def \Fix {{\text{Fix}\, }}
\def \Int {{\text{Int}\, }}
\def \Det {{\text{Det}\, }}
\def \Ker {{\text{Ker}\, }}
\def \Imm {{\text{Im}\, }}
\def \Sing{{\text{Sing}\, }}
\def \Tr {{\text{Tr\,  }}}
\def \sgn{{\text{sgn}\,}}
\def\AA{{\Cal A}}
\def\BB{{\Cal B}}

\def\DD{{\Cal D}}

\def\FF{{\Cal F}}

\def\KK{{\Cal K}}
\def\LL{{\Cal L}}
\def\MM{{\Cal M}}
\def\NN{{\Cal N}}

\def\PP{{\Cal P}}
\def\QQ{{\Cal Q}}
\def\RR{{\Cal R}}
\def\SS{{\Cal S}}

\def\VV{{\Cal V}}
\def\XX{{\Cal X}}

\def\today {\ifcase\month\or January \or February \or March \or
April \or May \or June
\or July \or August \or September \or October \or November \or December
\fi
\number\day~\number\year}


\topmatter
\title\nofrills
Fonctions z\^eta dynamiques\\
(Dynamical Zeta Functions)
\endtitle

\author
Viviane Baladi (Cours de DEA, Orsay, 2002)
\endauthor

\date
Rough course notes \quad June 2002
\enddate

\address
V. Baladi: CNRS,
IHES, 35, route de Chartres, F-91440 Bures-sur-Yvette, FRANCE
\endaddress
\email
baladi\@ihes.fr
\endemail

\endtopmatter

\document
\head
1. Introduction
\endhead

Dynamical zeta functions or dynamical determinants
are power series $\zeta(z)$, respectively $d(z)$, which are
constructed from (weighted)  periodic orbit data arising
from a, say, discrete-time dynamical system $f:M\to M$
and a function $g :M \to \complex$,
and which play the part of a (generalised) Fredholm determinant
for the transfer operator $\LL$ (on a suitable Banach space)
associated to $f$ and the weight $g$,
in the sense that they define a meromorphic, respectively
holomorphic, function in some domain where their poles
(respectively zeroes) are in bijection with the inverse
eigenvalues of $\LL$ in this domain. Although we shall not
explain this here, the spectral properties of $\LL$ are often
closely related to the statistical properties of the dynamical
system.

\remark{Exercise 0}
Let $L$ be a finite matrix with complex coefficients. Check that
$$
\det (\Id-zL) = \exp- \sum_{n=1}^\infty {z^n \over n}
\Tr L^n \, .
$$
\endremark

\smallskip

In this introduction, we shall study the case of a one-dimensional
dynamical system (i.e. a transformation of a compact interval) and
see how far Exercise~0 can take us. The course will then be devoted
to a presentation of more sophisticated arguments, inspired initially
by the work of Milnor and Thurston, which will allow us to treat
completely the one-dimensional situation.  We shall
then discuss much more recent results in higher dimensions.

\medskip

\subhead 1.1 The one-dimensional setting -- the transfer operator $\LL$
\endsubhead

Let $I=[0,1]$ be the unit interval (one could take any other compact
interval) and let $f: I \to I$ be a continuous map which is piecewise
monotone and piecewise $C^1$ with inverse branches having a
derivative of bounded variation. This means that we assume
that there is a  partition of $I$ into $N$ nontrivial subintervals
$I_j=[a_{j}, a_{j+1}]$, $j=0, \ldots, N-1$ (here we only consider the
case of finite $N$) such that:
\roster
\item the restriction of $f$ to $I_j=[a_j, a_{j+1}]$ is strictly
monotone for each $j=0, \ldots, N-1$;
\item the restriction of $f$ to $I_j$ extends
to a strictly monotone $C^1$ map to a small neighbourhood of $I_j$,
and the absolute value of the derivative
of the inverse $\psi_j$ of this map
$$
g_j :=\chi_{f(I_j)} \cdot  \psi_j'=
\chi_{f(I_j)} \cdot {1\over |f' \circ f|_{I_j}^{-1}  }
$$
is of bounded
variation (on the closure of $f(I_j)$).
Here, we are slightly abusing notation and we have set
$f'(a_i) =\lim_{x < a_i, x\to a_i} f'(x)= f'(a_i-)$ for
$i=1, \ldots N$, and $f'(a_0) = f'(a_)+)$.
\endroster

\noindent We shall {\it not} assume that the intervals
$I_j$ are maximal for the monotonicity property~ \therosteritem{1}.

\smallskip
We recall for the convenience of the reader that a
function $g : \real \to \complex$ is of bounded variation, noted
$g \in BV$ if
$$
\var_\real g =\sup_{m, t_i}  \{ \sum_{i=0}^m |g(t_i)-g(t_{i+1})| \, ,
t_0 < t_1 < \ldots < t_{m+1} \} < \infty \, .
$$
If $J \subset \real$ then $g$ is of bounded variation on $J$,
noted $g \in BV(j)$ if $\var_J (g)<\infty$ where $\var_J(g)$
is the above supremum restricting the partitions to
$t_i \in  J$.

\remark{Exercise 1}
Check that assumption \therosteritem{1} implies that $\sup_j g_j$
is finite.
\endremark

\smallskip

Recall that $f_*(\mu)$, if $\mu$ is a finite complex Borel measure,
is defined by $f_*(\mu)(E)= \mu(f^{-1}(E))$ for all Borel sets $E$.
Exercise ~1 can be used to prove that the {\it transfer operator}
(also called Ruelle operator, or, in this specific context,
Perron-Frobenius or density transformer operator) defined
as an operator on $L^1=L^1(I,Leb)$ by
$$
\int \psi (\LL \varphi) \, dLeb = \int \varphi (\psi \circ f) \,
dLeb \, , \forall \varphi \in L^1\, , \psi \in L^\infty (I, Leb)\, ,
$$
or, equivalently
$$
(\LL\varphi )\, dLeb = f_*(\varphi \, dLeb) \, ,
$$
or finally
$$
\LL \varphi(x)=\sum_{fy=x} {\varphi(y) \over |f'(y)|}
=\sum_{j=0}^{N-1} g_j(x) \varphi \circ \psi_j(x) \, ,
$$
is bounded, and that
$$
| \int_I \LL\varphi \, dLeb| \le \int_I |\varphi| \, dLeb\, ,
$$
i.e. the norm of $\LL$ on the Banach space $L^1$ is
(at most) $1$.

\remark{Exercise 2}
Find an example of $f$ such that $\LL$ does not preserve the
Banach space $C^0(I)$ of continuous functions on $I$.
\endremark

\remark{Exercise 3}
Show that if there is  a nonnegative $\varphi_0 \in L^1$
with $\LL\varphi_0 = \varphi_0$  (and $\int \varphi_0 \, dLeb> 0$)
then the measure $\mu_0=\varphi_0 \, dLeb/\int \varphi_0\, dLeb$
is an (absolutely continuous) $f$-invariant
(probability) measure, i.e. $f_*(\mu_0)=\mu_0$.

Absolutely continuous invariant measures $\mu_0$ are especially
interesting when they are ergodic. Indeed, the Birkhoff ergodic
theorem then says that  for a set of {\it positive
Lebesgue measure} of initial conditions $x$ we have
($\delta_y$ denotes the Dirac mass at $y$ and the convergence is
in the weak* topology):
$$
\lim_{n\to \infty}{1\over n}
\sum_{k=0}^{n-1} \delta_{f^k x} = \mu_0 \, .
$$
This is often interpreted as an indication that such a measure
$\mu_0$ is a natural, or physical, measure for $f$.
\endremark

\medskip
We shall not discuss this here, but one can prove for example
if $\sup_j \inf g_j < 1$ that:
\roster
\item
$1$ is indeed an eigenvalue (of finite multiplicity) of $L^1$,
and the spectrum of $L^1$ on the unit circle consists in roots
of unity $e^{2 i  k \pi/k_0}$, $k=0, \ldots, k_0-1$ for some
integer $k_0\ge 1$; in fact, $1$ is a simple eigenvalue
if and only if $\mu_0$ is ergodic and in this case $k_0=1$ if and only
if $\mu_0$ is mixing;
\item
every point in the open unit disc is an eigenvalue of infinite
multiplicity of $\LL$ acting on $L^1$
(the spectrum of $\LL$ is therefore the entire closed unit
disc).
\endroster

It appears that $L^1$ is not very suitable to obtain
spectral information reflecting finer statistical properties
(stability of the absolutely continuous invariant
measure under small  deterministic or
probabilistic perturbations, exponential decay of correlations
for suitable observables, central limit theorem, etc.) of the
dynamics. In some sense, $L^1$ is too big a Banach space and
we should find a smaller invariant Banach space. (Note that
the Hilbert space $L^2$ suffers from the same ``problems'' as
$L^1$.) In our one-dimensional framework, the most natural candidate is
the Banach space $BV=BV(I)$ of functions of bounded variation
on $I$. (If we had assumed that the partition satisfies a Markov
property -- see below -- it would be possible to consider other
choices.)

\medskip

\subhead 1.2 The transfer operator acting on $BV$: quasicompactness
\endsubhead

\smallskip
Let us consider the Banach space
$$
BV=BV(I) = \{ \varphi : I \to \complex\, , \var_I \varphi < \infty\} \, ,
$$
endowed with the norm
$$
\|\varphi\|_{BV(I)} = \var_I \varphi + \sup_I |\varphi|  \, .
$$
(The supremum term is here to distinguish constant functions on
$I$, it could be replaced e.g. by the $L^1$ norm, or also by
substituting $\var_I$ by $\var_\real$.)

In fact, it will be more convenient to consider the quotient
$$
\BB = BV/\NN \, ,
 \|\varphi\|=\|\varphi\|_\BB=\inf\{ \|\phi\|_{BV(I)}\mid \phi-\varphi \in \NN \}\, ,
$$
where $\NN$ is the space of complex-valued functions on $I$
which vanish except on an
at most countable set.

\remark{Exercise 4}
Show that $BV(I)$  and $\BB$ are indeed Banach spaces (i.e. they
are complete) and that $\LL$ maps $BV(I)$ into $BV(I)$
boundedly.
\endremark

\smallskip

We are now ready to state and prove our first result:

\proclaim{Theorem 1 (Quasicompactness of $\LL$ on $\BB$)}
Let $f$  and $\LL$ be as above.  Then $\LL$ is a bounded
operator on $\BB=BV(I)/\NN$, and, outside of the closed disc
of radius
$$
\widehat R := \limsup_{n \to \infty} \sup_x
\left ( {1\over |(f^n)'(x)|} \right )^{1/n}\, ,
$$
the spectrum of $\LL$ consists in isolated eigenvalues of
finite multiplicity.
\endproclaim

One can prove additionally that the spectrum of $f$ on $BV(I)$
and on $\BB$ coincide outside of the disc of radius $\widehat R$,
we shall not do this here.

\demo{Proof of theorem 1}
Since the transfer operator is a sum of
operators of composition and multiplication,
and since the variation is essentially the $L^1$
norm of the distributional derivative,
the ingredients are the Leibniz formula (derivative of a product)
and the chain rule (derivative of a composition).
We present a conceptual proof due to Ruelle,
the starting point of which is to replace elements of
$\BB$ by Radon measures:
\enddemo

\proclaim{Lemma 0 (Bounded variation and Radon measures)}
The Banach spaces
$$
\BB'=\{ \varphi \in \BB \mid \varphi(0+)=0\}
$$
and
$$
C^0(I)^* = \{ \nu : C^0(I)\to \complex\mid \hbox{ linear and continuous} \}
$$
are isomorphic, the Banach space isomorphism being
given by the distributional derivative (Stieltjes measure
associated to a function of bounded variation)
$$
d \varphi (c,d] = \varphi (d+)-\varphi(c+)\, , (c,d] \subset I\, , \quad
\varphi \in BV \, .
$$
The inverse of $d$ will be denoted by $\SS$, and satisfies
$\SS \mu(x)=\mu([0,x])$.
\endproclaim

Recall that if $\nu$ is a Radon measure and $\varphi$ a
bounded function the Radon measure $\varphi \nu$ is
defined by $\varphi \nu(\psi) = \nu(\varphi \psi)$ where
we use that $C^0(I)^*$ is the space of bounded complex Borel
measures on the compact metric space $I$.

\proclaim{Lemma 1 (Leibniz formula in BV/integration by parts)}
Let $\varphi_1$, $\varphi_2 \in BV(I)$. In the quotient
$\BB=BV(I)/\NN$ we may take representatives which
are continuous at $a_0=0$ and which only have regular
discontinuities $2 \varphi_i(x)=\varphi_i(x+)+\varphi_i(x-)$.
Then, for these representatives
$$
d(\varphi_1 \varphi_2) = \varphi_1 d(\varphi_2) + \varphi_2 d\varphi_1\, ,
$$
in the sense of Radon measures.
\endproclaim

\proclaim{Lemma 2 (Change of variables)}
Let $J \subset I$ be an interval, $\psi : J \to \psi(J)$ be
a homeomorphism, and $\varphi \in BV(I)$. Then
$$
\chi_J \, d(\varphi \circ \psi) =
(\epsilon \chi_J ) (\psi^{-1})_* (d\varphi)\, ,
$$
where $\epsilon = +1$ if $\psi$ preserves the orientation
and $-1$ if $\psi$ reverses the orientation.
\endproclaim

The sign in Lemma 2 comes from the fact that if $(c,d] \subset J$, e.g.
then
$$
d(\varphi \circ J) (c,d] = \varphi (\psi(d+)) - \varphi(\psi(c+))
$$
while
$$
(\psi^{-1})_* (d\varphi) (c,d] = (d\varphi) \psi(c,d] =
\cases
(d\varphi) (\psi c,\psi d]
= \varphi(\psi(d+))- \varphi(\psi(c+))& \hbox{if}\, \epsilon > 0 \, ,\cr
(d\varphi) [\psi d,\psi c)
=\varphi(\psi(c+))- \varphi(\psi(d+))& \hbox{if} \,\epsilon < 0 \, .
\endcases
$$

The complete proofs of Lemmas 0--2 are to be found, e.g.,
in [DS1] or in [Ba2].

\demo{Back to the proof of Theorem 1}
Our first step is to replace the transfer operator $\LL$ on $\BB$
by $\MM$, the rank-one perturbation of $\LL$ given by
$$
\MM\varphi= \LL \varphi - \LL\varphi(0+) \, .
$$
The arguments below will show that if we can prove the
claim on $\MM$ acting on $\BB$, then it will be automatically
satisfied for $\LL$ acting on $\BB$.
The next observation is that $\MM$ maps $\BB$ into
$\BB'$, so that it suffices to analyse the spectrum of
$\MM$ on $\BB'$ (recall the definition of the nonzero  spectrum
and write $(\Id-z\MM)^{-1} = (\Id -z \MM)^{-1} z \MM +\Id$).

The operator on $C^0(I)^*$ conjugated to $\MM:\BB'
\to \BB'$ by the isomorphisms of Lemma~0 is
$d \circ \MM \circ \SS$. Applying Lemmas~1 and 2,
it is not difficult to see that $d\MM \SS$ can be decomposed
as
$$
d\MM\SS = \widehat \MM + \widehat \NN \SS \, ,
$$
where, setting $\epsilon_i=+1$ if $\psi_i$ is increasing
and
$\epsilon_i = -1$  if $\psi_i$ is decreasing
$$
\eqalign{
\widehat \MM (\mu) (\varphi)
&= \sum_{i=0}^{N-1}
\chi_{(a_i, a_{i+1}]} \epsilon_i ( ( g_i \varphi) \circ \psi_i^{-1} )
\, d\mu \cr
&=\int \epsilon_f {\varphi \circ f \over |f'| } \, d\mu\, ,
}
$$
where $\epsilon_f = \epsilon_i $ on $(a_i, a_i+1)$
and $0$ on the $a_i$s,
or, introducing the more compact notation $g |_{(a_i,
a_{i+1}]}:= g_i \circ \psi_i^{-1}$ (and $g(a_0)=g(a_0+)$),
$$
\widehat \MM (\mu )= (\epsilon_f \cdot g) f_*(\mu) \, ,\quad
\widehat \MM : C^0(I)^* \to C^0(I)^* \, ,
$$
and
$$
\widehat \NN( \varphi)=
\sum_{i=0}^{N-1} d(g_i) \varphi \circ \psi_i \, ,\quad
\widehat \NN : \BB' \to C^0(I)^* \, .
$$
(We shall not need this but
note that $\widehat \NN$ can be viewed as $d\widehat \MM_0-
\widehat \MM d$ where $\widehat \MM_0$ is $\widehat \MM$
``acting on functions.'' )

It is not difficult to check that the spectral radius of $\widehat \MM$
on $C^0(I)^*$ is not larger than $\widehat \RR$.
The next important step is encapsulated in a last sublemma:
\enddemo

\proclaim{Lemma 3}
Let $\mu_0$ be a fixed Radon measure. Then the linear
operator on $C^0(I)^*$
$$
\mu \mapsto (\SS\mu) \cdot \mu_0
$$
is compact
\endproclaim

\demo{Sketch of the proof of Lemma 3}
Since $\mu \mapsto \SS \mu$ is bounded from $C^0(I)^*$
to $\BB'$ it suffices to see that the operator $\KK$ from
$\BB'$ to $C^0(I)^*$ defined by
$$
\KK(\varphi)= \varphi \cdot \mu_0
$$
is compact. For this, we need to prove that any sequence
$\KK(\varphi_n)$ with $\var \varphi_n\le 1$ admits
a convergent subsequence. This can be deduced from the fact
that $\KK$ can be approached by a sequence of finite-rank
operators. To construct this sequence, note that for any
$\delta > 0$ there is a  finite
partition $c_0 < \cdots < c_i < \cdots < c_M$
of $I$ such that the norm of the Radon
measure $\chi_{(c_i, c_{i+1})} \mu_0$
is smaller than $\delta$ for each $0\le i< M$.
Setting
$$\Pi_\delta \varphi|_{(c_i, c_{i+1})}= \varphi(c_i+)$$
for all $0 \le i < M$, it is clear that $\Pi_\delta:\BB'\to \BB'$
is a finite-rank (rank $M$, in fact, projection), so that
the operator $\KK_\delta := \KK \Pi_\delta$ is
also finite-rank from $\BB'$ to $C^0(I)^*$.
Finally, it is not difficult to show that
$$
\|(\KK \varphi - \KK_\delta \varphi)\|=
\|\varphi  \mu_0 - \Pi_\delta \varphi \mu_0\|
\le \delta \var_I \varphi \, .
$$
(A full proof can be found, e.g., in  [Ba2].)
\enddemo

\demo{End of the proof of Theorem 1}
Applying Lemma~3 to the finite set of measures of
the form $\mu_0=d(g_i)$, we see that $\widehat \NN \SS$
is a compact operator. If the reader knows that a perturbation
of a bounded operator of spectral radius $\rho$ by a compact operator
can only add isolated eigenvalues of finite multiplicity
outside of the disc of radius $\rho$ (``compact perturbations
do not change the essential spectral radius''), he or she should
be satisfied. Otherwise, let us proceed with the proof,
decomposing the resolvent as
$$
(\lambda \Id - d\MM \SS)^{-1}
= (\lambda \Id -( \widehat \MM + \widehat \NN \SS))^{-1}
=(\lambda \Id -\widehat \MM)^{-1}
(\Id-\widehat \NN \SS (\lambda\Id - \widehat \MM)^{-1})^{-1}\, .
$$
If $|\lambda| > \widehat R$, the resolvent $(\lambda \Id -\widehat \MM)^{-1}$
of $\widehat \MM$
is a bounded operator (depending holomorphically
on $\lambda$ in the domain $C_{\widehat \RR}=\{ |\lambda| > \widehat R\}$). Therefore, the operator
$$
\QQ(\lambda)=
\widehat \NN \SS (\lambda\Id - \widehat \MM)^{-1}\, ,
$$
being the composition of a bounded operator and a compact
operator, is compact. It also depends holomorphically on
$\lambda$ in $C_{\widehat \RR}$.
Our aim is therefore to study the set
$$
\Sing=\{ \lambda\in C_{\widehat \RR} \mid
1 \hbox{ is an eigenvalue of } \QQ(\lambda) \} \, .
$$
Our first remark is that $\Sing$ is a discrete subset
of (a bounded subset of) the complex plane.
Indeed, as $|\lambda|\to \infty$, the spectral radius
of $\QQ(\lambda)$ goes to zero, so that $\QQ(\lambda)$
cannot have an eigenvalue $1$ for all
$\lambda$  in the connected domain $C_{\widehat R}$.
Since the
nonzero eigenvalues of a family of compact operators
depending analytically on a parameter $\lambda$
are either constant or take any fixed value on a
discrete set, we are done. (The last claim is analogous
to the corresponding result for finite matrices --
simple eigenvalues depend analytically on
analytic perturbations, multiple eigenvalues
can have at worse algebraic (roots) singularities --- and both
statements can be found in [Ka,
II.\S1 and VII.\S1], e.g.)

It remains to be seen that any point in the discrete set
$\Sing$ is an eigenvalue of finite multiplicity of
$\widehat \MM + \widehat \NN \SS$.
If $\lambda \in \Sing$,
it is not difficult to associate to the fixed function
$\varphi_\lambda$ of $\QQ(\lambda)$  an
eigenfunction of $\widehat \MM + \widehat \NN \SS$
for the eigenvalue $\lambda$. (This is left
as an exercise to the reader.) This does not completely
end our task, since this eigenvalue could in principle
have infinite multiplicity. In order to finish the proof,
we present a few reminders about the theory
of spectral projectors associated to isolated points in
the spectrum of an operator (see, e.g., [Ka]).

So let $\lambda_0$ be an isolated point in the spectrum of
a bounded linear operator $L$ on a Banach
space (we are not assuming that
$\lambda_0$ is an eigenvalue). This implies that there is
a nontrivial complex disc $D_\gamma(\lambda_0)$
centered at $\lambda_0$ which does
not intersect any other point in the spectrum of $L$.
Letting $\gamma=\gamma(\lambda_0)$ be the path corresponding to going
along once the circle
bounding this disc counterclockwise, we define
a bounded operator on our Banach space by
$$
P^L_{\lambda_0}=
{1\over 2i\pi} \oint_\gamma (\lambda \Id -L)^{-1} \, d\lambda \, .
$$
(To check that the sign is correct, consider $L\equiv 0$ and $\lambda_0=0$
and verify that $P^L_0=\Id$.)
Let us verify that $P^L_{\lambda_0}$ is a
projector, i.e. $(P^L_{\lambda_0})^2=P^L_{\lambda_0}$.
For this, the first remark is that if $\gamma'$ is the circle
centered at $\lambda_0$ and of radius one-half the radius
of $\gamma$ (for example) then since
$(\lambda \Id - L)^{-1}$ is holomorphic in the annulus
bounded by the two circles we can also write
$$
P^L_{\lambda_0}=
{1\over 2i\pi} \oint_{\gamma'} (\lambda' \Id -L)^{-1} \, d\lambda' \, .
$$
Therefore, using the easily checked ``resolvent identity''
$$
(\lambda \Id-L)^{-1}- (\lambda'\Id-L)^{-1}
= (\lambda-\lambda') (\lambda\Id-L)^{-1} (\lambda'\Id-L)^{-1} \, ,
$$
we find
$$
P^L_{\lambda_0}P^L_{\lambda_0}
= {1\over (2i\pi)^2}
\int_\gamma \int_{\gamma'} (\lambda-\lambda')^{-1}
\left [(\lambda\Id -L)^{-1} -(\lambda'\Id -L)^{-1} \right ]
\, d\lambda'\, d\lambda \, .
$$
We finish by observing that
$$
{1\over 2i\pi} \int_{\gamma'} (\lambda-\lambda')^{-1} \, d\lambda'=0\, ,
$$
and
$$
{1\over 2i\pi}
\int_\gamma  (\lambda-\lambda')^{-1} \, d\lambda=1\, .
$$

$P^L_{\lambda_0}$ being a projector,
it follows that $\Id -P^L_{\lambda_0}$ is also a
projector, and clearly the two projectors are orthogonal
(i.e $P^L_{\lambda_0}(\Id -P^L_{\lambda_0})=0=
(\Id -P^L_{\lambda_0})P^L_{\lambda_0}$).
Also, one easily checks that the definition implies
$LP^L_{\lambda_0}=P^L_{\lambda_0} L$, and, of course,
similarly for the other projector. Finally, one can show
that the nonzero spectrum of $L P^L_{\lambda_0} $ consists
in the single point $\lambda_0$, while the spectrum
of  $L (\Id -P^L_{\lambda_0} )$ does not intersect the
closed disc $D_\gamma$ centered at $\lambda_0$.
Now, in the case when $P^L_{\lambda_0}$ is finite
rank,  the operator $L P^L_{\lambda_0} $
acting on the finite-dimensional Banach space
$\Imm P^L_{\lambda_0}$ is of course finite rank, so
that its spectrum (which we already know is $\{\lambda_0\}$)
must be an eigenvalue of finite multiplicity.
(Note that $\Imm P^L_{\lambda_0}$ is the generalised
eigenspace, i.e. it does not always contain only  eigenfunctions
but also generalised eigenfunctions $\varphi$ such that
$(\lambda_0 \Id -L)^k \varphi =0$ for some $k\ge 2$
but $\ne 0$ for $k=1$. In particular, the dimension of
$\Imm P^L_{\lambda_0}$ is the algebraic multiplicity of
$\lambda_0$.)

Let us return to our specific problem, i.e. showing that
a point $\lambda_0$, such that $1=\rho(\lambda_0)$ is an eigenvalue of
$\QQ(\lambda_0)$, is an eigenvalue of finite multiplicity
of $d\MM \SS$.  Since the nonzero spectrum of a compact
operator consists in isolated points, and since the  small
perturbations $\QQ(\lambda)$
in operator norm of $\QQ(\lambda_0)$ produce small perturbations
$\rho(\lambda)$ of our isolated eigenvalue $1$, up to taking
a  smaller isolating disc $D_\gamma(\lambda_0)$ for
$d\MM\SS$ and the spectral point $\lambda_0$,
we may find a disc centered at $1$, bounded by a
curve $\Gamma=\Gamma(\lambda_0, \gamma)$, such that for each
$\lambda$ in $D(\lambda_0)$, the
curve $\Gamma(\lambda_0, \gamma)$ does not intersect
the spectrum of $\QQ(\lambda)$.
In particular, the spectral projectors
$$P^{Q(\lambda)}_{\rho(\lambda)}
={1\over 2i\pi} \int_\Gamma (\rho\Id  - \QQ(\lambda))^{-1}\, d\rho \, ,
$$
will  have constant rank equal to the multiplicity of
$1$ for $\QQ_{\lambda_0}$ for all
$\lambda \in \gamma =\gamma(\lambda_0)$.

It will follow that the spectral projector
$$
P^{d\MM\SS}_{\lambda_0}=
{1\over 2i\pi} \int_\gamma (\lambda \Id - d\MM \SS)^{-1} \, d\lambda \, ,
$$
can be written as a path integral of a finite-rank operator.
Indeed, we may use the finite rank spectral projectors
$P^{\QQ(\lambda)}_{\rho(\lambda)}$ associated to the
perturbation of the eigenvalue $1$ for $\QQ_{\lambda_0}$ to refine
our previous decomposition of the resolvent of $d\MM\SS$ as:
$$
\eqalign{
(\lambda \Id - d\MM \SS)^{-1}
&=(\lambda \Id -\widehat \MM)^{-1}
(\Id -\widehat \NN \SS (\lambda\Id - \widehat \MM)^{-1})^{-1}\cr
&=(\lambda \Id -\widehat \MM)^{-1}
(\Id-\QQ(\lambda) P^{Q(\lambda)}_{\rho(\lambda)})^{-1}
P^{Q(\lambda)}_{\rho(\lambda)}\cr
&\qquad\qquad
+(\lambda \Id -\widehat \MM)^{-1}
(\Id-\QQ(\lambda))^{-1}(\Id -P^{Q(\lambda)}_{\rho(\lambda)})\, .
}
$$
The second term in the above decomposition being holomorphic in
the disc bounded by $\gamma$, the corresponding path integral
vanishes. The first term is the composition of a bounded
operator and a finite-rank operator, it is thus finite-rank.
It also depends holomorphically on $\lambda$ on
$\gamma$ (and meromorphically on $\lambda$ in the disc bounded by $\gamma$).

To finish, it suffices to note that the path integral of finite-rank
operators being compact, the projector $ P^{d\MM\SS}_{\lambda_0}$
is compact and therefore necessarily finite-rank.
\qed
\enddemo

\remark{Exercise 6}
Note that since the dual of $\LL$ acting on the dual of
$\BB$ preserves Lebesgue measure, the operator $\LL^*$ has
a fixed point in $\BB^*$. Assume that $\widehat R < 1$.
Since the spectrum of $\LL$
on $\BB$ outside of the disc of radius $\widehat R$ consists
only of isolated eigenvalues of finite multiplicity, show that
any eigenvalue
of $\LL^*$ in this domain must be an eigenvalue of $\LL$.
It follows that $\LL$ has a fixed
point in $\BB$.
\endremark

\medskip

\subhead 1.3 The dynamical zeta function in the Markov expanding affine case
\endsubhead

\smallskip

We shall now combine Theorem~1 and Exercise~0 to obtain
a result on the dynamical zeta function of $f$, but only under three
additional assumptions which are quite restrictive. In some
sense, the purpose of the course is to show how one can get rid
of these assumptions. Here  are the first two:

\roster
\item
We suppose that $f'|_{(a_i, a_{i+1})}$ is constant for each $i$;
\item
and that the partition of $I$ into $N$ intervals $I_j$
satisfies the Markov property: if $f(I_j)\cap \Int I_\ell
\ne \emptyset$ then $\Int I_\ell \subset f(I_j)$.
(In other words, each $\overline {f(I_j)}=$ can be written exactly
as a union of $\overline{I_\ell}$s.)
\endroster

Let us consider the $N$-dimensional vector subspace
of $\BB$ defined by
$$
\VV =\{ \varphi \in \BB \mid \varphi|_{(a_i, a_{i+1})}
\hbox{ is constant } \} \, .
$$
It is not difficult to check that $\LL$ maps $\VV$ into itself.
Also, introducing the $N \times N$
Markov (or transition) matrix associated to $f$:
$$
A_{jk} =  1 \hbox{ if } f(I_j) \cap \Int I_k \ne \emptyset \, ,
\quad
A_{jk} =  0 \hbox{ if } f(I_j) \cap \Int I_k = \emptyset \, ,
$$
it is easy to verify that the matrix of $\LL|_{\VV}$ in the
standard basis is given by the matrix $A_g$ defined by
$$
(A_g)_{jk }=A_{jk} g_j \, .
$$
The spectrum of $\LL$ on $\VV$ is thus a well defined set
of eigenvalues (of finite multiplicity). Note that we do not claim
that this sets intersects the complement of the disc of
radius $\widehat R$ (but see Exercise~ 6 and Lemma~ 4).

Our first result is:

\proclaim{Lemma 4}
Assume \therosteritem{1--2}.
Outside of the closed disc of radius $\widehat R$, the
spectrum of $\LL$ on $\BB$ coincides with the spectrum of
$\LL$ on $\VV$.
\endproclaim

\demo{Proof of Lemma 4}
Since $\VV \subset \BB$,
if $\varphi \in \VV$ is an eigenfunction for $\lambda$ and $\LL$ then
$\lambda$ is also an eigenvalue for $\LL$ acting on $\BB$.
(Note that this is also true if $|\lambda|\le \widehat R$.)
To show the reverse inclusion, let us suppose
that there is $\lambda$ with $|\lambda|>\widehat R$
and a nonzero $\varphi \in \BB$ with $\LL \varphi = \lambda \varphi$
(we know by Theorem~1 that in this domain the spectrum of
$\LL$ on $\BB$ consists
in eigenvalues). Then,  using the operator $\MM$ from
the proof of Theorem~1, we have for all $k \in \natural$
$$
\varphi - {\varphi(0+)\over \lambda}
= \lambda^{-k} \MM^k \varphi \, .
$$
Let us rewrite the right-hand-side
of the above equality in the coordinates given by Lemma~0,
using the notation there
$$
\lambda^{-k}  \widehat \MM^k d \varphi +
\lambda^{-k}
\sum_{m_i \, n_i \in \{0, 1\}  \, , \sum m_i+n_i=k\, \sum n_i \ge 1}
\prod_{i=1}^k \widehat \MM^{m_i} (\widehat \NN \SS)^{n_i} (d \varphi) \, .
$$
For any $\widetilde R > \widehat R$ there is $C > 0$
so that  $\| \widehat \MM ^k(d\varphi)\|
\le C \widetilde R^k \|\varphi\|$, so that, taking
$\widetilde R > |\lambda| >\widehat R$, we see that the
first term in the above decomposition goes to zero as $k$ goes
to infinity. Let us then concentrate on
the second term.
Our Markov and piecewise affine assumptions imply that
each of the $N$ measures $d(g_i)$
is a (positive) linear combination
of the two Dirac masses at $f(a_i)=a_{u(i)}$ and $f(a_{i+1})=
a_{v(i)}$ (by construction, the weight of $\delta_{a_0}$
vanishes).
Multiplying $d(g_i)$ by an arbitrary element $\psi$ of $\BB$
amounts to changing  the coefficients of this linear
combination (using $\psi(a_{u(i)}\pm)$ and $\psi(a_{v(i)}\pm)$).
Therefore, it suffices to analyse the action of
$\widehat \MM^\ell$ ($\ell \ge 1$) on a linear combination
of Dirac masses at the $a_j$s:
Using again the Markov assumption, we see that
$\widehat \MM^\ell$ is still a linear combination at the
endpoints $a_j$.
Putting everything together, we see that
$d(\varphi-\varphi(0+)/\lambda)$ must be a linear combination of Dirac masses
at the $a_j$s, so that $\varphi$ belongs to $\VV$, as claimed.
\qed
\enddemo

From now on until the end of the introduction, we make our third and final
additional assumption:

(3) $\widehat R < 1$.

\medskip

We shall use the notation $\Fix f$ for the set
$$
\Fix f=\{ x\in I \mid f(x) = x\}\, ,
$$
and similarly for $\Fix f^n$ for each nonnegative integer.

\remark{Exercise 7}
Prove that assumption \therosteritem{3} implies that
for each $n$ the set $\Fix f^n$ is finite. (We shall obtain more precise
information on this set soon.)
\endremark

Exercise 7 allows us to define the weighted dynamical zeta function
of $f$ by the following formal power series (recall the notation
$g$ from the proof of Theorem~1):
$$
\zeta(z)=\zeta_{f,g}(z)
=\exp \sum_{n =1}^\infty
{z^n \over n} \sum_{x \in \Fix f^n} \prod_{k=0}^{n-1} g(f^k (x))\, .
$$

Here is the final result of the introduction:

\proclaim{Theorem 2 (Zeta function and spectrum in the Markov
expanding affine case)}
Assume \therosteritem{1--3} from Subsection~ 1.1 and
\therosteritem{1--3} from Subsection~ 1.3.
Then the weighted zeta function $\zeta(z)$ is meromorphic in the open
disc of radius $\widehat R^{-1}>1$, where its poles are exactly the inverse
eigenvalues (outside of the closed disc of radius
$\widehat R$) of $\LL$ acting on $\BB$.
The order of each such pole $z$ coincides with the algebraic
multiplicity of the
corresponding eigenvalue $1/z$.
\endproclaim

\smallskip
\demo{Proof of Theorem 4}
The Markov assumption allows us to construct a symbolic model for $f:I\to I$
which is a subshift of finite type (SFT). Let us define this SFT, recalling
the transition matrix $A$. Consider the
set $\Sigma$ of one-sided sequences with coefficients in
the finite alphabet $\{0, \ldots, N-1\}$ ,
and the subset $\Sigma_A\subset\Sigma$ of sequences with $A$-admissible
transitions, i.e.
$$
\Sigma_A=\{ t \in \Sigma \mid A_{t_i t_{i+1}}=1\, ,  \forall i \in \integer \}\, .
$$
(This is a compact set for the product topology arising from the
discrete topology on our finite alphabet.)
The one-sided shift to the left $(\sigma(t))_i = t_{i+1}$ leaves $\Sigma_A$
invariant. We next construct a semi-conjugacy between $\sigma|_{\Sigma_A}$
on $\Sigma_A$ and $f$ on $I$, i.e. a surjective map $\pi : \Sigma_A\to I$
with $f \circ \pi = \pi \circ \sigma$ on $\Sigma_A$. For this one first
observes that for each $t \in \Sigma_A$ the set
$
\cap_{i=0}^\infty f^{-i} I_{t_i}
$
is a single point in $I$. (This set is nonempty because the sequence
is admissible, and it has zero diameter because we assumed $\widehat R<1$.)
Setting
$$
\pi(t) = \cap_{i=0}^\infty f^{-i} I_{t_i}\, ,
$$
we obtained the desired conjugacy. To check surjectiveness, note that
each trajectory of each point $x \in I$, i.e. each admissible sequence of symbols
$t_i$ with $f^ i (x) \in\overline{ I_{t_i}}$ gives $t\in \Sigma_A$ with
$\pi(t)=x$. There may be an ambiguity if $f^i (x) = a_j$  for some $i, j$,
so that the map $\pi$ is not injective in general.

In order to obtain a conjugacy (i.e. a bijection making the diagram commute),
it is convenient to slightly modify the original interval map $f$ on $I$ by
``doubling'' the $N-1$ points $a_1, \ldots, a_{N-1}$ and all their countably
many preimages $f^{-k}a_i$,
$k \ge 1$. Between each such pair of doubled points we introduce a
small interval of length, say, $\epsilon/(N^k k^2)$, in such a way
that the total added length is finite. This allows us to embed our Cantor set
$\widehat I$ into a compact
interval of the real line. Abusing slightly notation,
$I \subset \widehat I$ and the closed intervals $\overline I_j$ are disjoint in $\widehat I$.
We may extend $f$ to $\widehat I$, the new map $\hat f$ being just $f$ in the
interior of  each $I_j$s and being set to
$f(a_i +)$,  respectively $f(a_i -)$ in the new right or left boundaries.
Similarly, we extend the weight $g$ to $\widehat I$ by taking
the appropriate left or right limit. It should be clear that $\pi$ is now a bijection between
the Cantor sets
$\Sigma_A$ and $\widehat I$, such that $\hat f \pi=\pi \sigma$ on $\Sigma_A$.

Let us next analyse the weighted zeta function of $\sigma$ and the weight
$\hat g \circ \pi$, i.e.
$$
\zeta_{\sigma,\hat g \pi} (z)=
\exp\sum_{n=1}^\infty{z^n \over n}\sum_{t\in \Fix \sigma^n }\prod_{k=0}^{n-1}
\hat g(\pi(\sigma^k(t))\,  .
$$
We claim that
$$
\sum_{t\in \Fix \sigma^n }\prod_{k=0}^{n-1}
\hat g(\pi(\sigma^k(t))=\Tr A_g^n \, .
$$
Indeed, the above equality is obvious for $n=1$. More generally we have
$$
(A_g^n)_{i i}
=\sum_{i=i_0, \ldots, i_{n-1}\, ,(\vec \imath)^\infty \in \Sigma_A}
\prod_{k=0}^{n-1} g_{i_k}\, ,
$$
where $(\vec \imath)^\infty$ means the finite
(length-$n$) sequence $\vec \imath$ repeated infinitely many times.
All fixed points of $\sigma^n$ are obtained that way, and we have
$g(\pi(\sigma^k(\vec \imath^\infty))=g_{i_k}$.

Applying Exercise 0 to the finite matrix $A_g$, it follows that
$$
\zeta_\sigma (z)=(\det(\Id-zA_g ))^{-1} \, .
$$
By Lemma 4, we know that the eigenvalues of $A_g$ outside of the
disc of radius $\widehat R$ are in bijection with the spectrum of $\LL$,
outside the disc of radius $\widehat R$, acting on $\BB$.

By construction
$$
\zeta_{\hat f, \hat g}  (z)=\zeta_{\sigma,\hat g \pi} (z)\, ,
$$
so that, it suffices to understand the relation between $\zeta_{f,g}(z)$
and $\zeta_{\hat f,\hat g}(z)$ to end the
proof. More precisely it suffices to show that their
ratio is a nonzero holomorphic function in the  disc of radius $\widehat R^{-1}$.
Clearly, the periodic points of $f$ whose orbits do not meet any of the
$a_i$s are in bijection with the periodic points of $\hat f$ whose orbits do not
meet $a_0$, $a_N$, and any of the twins $a_{i\pm}$, $i=1, \ldots, N-1$;
also the contributions of these ``good'' periodic points to the respective zeta functions coincide.
(Recall the Markov assumption.)
So let us consider one of the finitely many possible
periodic points $a_i =f^{p(i)} (a_i)$ (assuming that $p(i)\ge 1$
is minimal for the fixed point property). There are three possibilities:
if $a_i$ is a local extremum for  $f^{p(i)}$
then either $a_i +$  (if the extremum is a minimum) or
$a_i -$ (maximum), but not both, will be a periodic point
for $\hat f$, with minimal period equal to $p(i)$.
If $f^{p(i)}$ is increasing in a neighbourhood of $a_i$
(or if $i=0$), then both $a_i +$ and $a_i -$ are periodic of
minimal period $p(i)$ for $\hat f$. If $f^{p(i)}$ is decreasing in a
neighbourhood of $a_i$, then both $a_i+$ and $a_i-$ will be periodic
points for $\hat f$, but their minimal period will be
$2p(i)$. The analysis just made also describes the (finitely many)
periodic points of $\hat f$ whose orbits
meet $a_0$, $a_N$, or any of the twins $a_{i\pm}$, $i=1, \ldots, N-1$.

Let us consider one of these finitely many ``bad'' periodic orbits
$x=f^p x$, or  $\hat x=\hat f^{\hat p} \hat x$,
with $p, \hat p\ge 1$ its minimal period
and $\lambda=\prod_{k=0}^{p-1} g(f^k x)$
respectively $\hat \lambda=\prod_{k=0}^{\hat p-1} \hat g(\hat f^k \hat x)$
the associated weight. Clearly, the corresponding contribution to the weighted
zeta function is:
$$
\exp \sum_{m=1}^\infty {z^{mp}\over mp} \lambda^{m}=
(1-z^p \lambda )^{-1/p}\, , \hbox{ or }
\exp \sum_{m=1}^\infty {z^{m\hat p}\over m \hat p} \hat \lambda^{m}=
(1-z^{\hat p} \hat \lambda )^{-1/\hat p} \, .
$$
To finish, it suffices to observe that $\lambda \ge \widehat R^p$,
$\hat \lambda \ge \widehat R^{\hat p}$.
(This follows from the definition for $\lambda$, while a
short argument is required for $\hat \lambda$.)
\qed
\enddemo

\smallskip
\head 2. Kneading theory in dimension one
\endhead

\smallskip
\subhead 2.1 Introduction
\endsubhead

Although the ideas at the basis of the kneading theory
in these notes are present in the very classical  paper
of Milnor and Thurston [MT], which was written in the seventies,
they were only applied to weighted zeta functions and the analysis
of the spectra of transfer operators in the nineties. Before
that, other methods had been developed (in dimensions one and
higher, and under various assumptions of expansion, hyperbolicity,
and/or regularity) in the continuation of the Markov approach
of Section ~1. In this introductory section, we first give a very
brief and incomplete presentation of some of the results obtained
by these older methods between 1976 and now, referring to [Ba1, Ba3]
for more general surveys; we then give a very brief presentation
of the key result of Milnor and Thurston which inspired  the
new kneading approach.
\smallskip

{\bf The Markov appproach for piecewise monotone maps}

This approach consists in viewing the situation of the subshift of finite
type and locally constant weight as the paradigm, and trying to make
more general weighted dynamical systems fit into this model.

The first generalisation of Theorem~2 in Section~1 consists in
maintaining the assumption that the piecewise monotone interval
map $f$ is piecewise affine (or more generally consider a locally
constant weight $g_i$, which could be unrelated to the derivative),
but relaxing the Markov assumption. (Historically, the case of
Markov maps with non locally constant weights was studied earlier
by Mayer e.g.,
since the important Gauss map $x \mapsto \{ 1/x\}$ fits
in this framework -- see also the discussion
about analytic systems below.) A helpful tool here is the
Markov extension devised by Hofbauer in the seventies. This
associates to any piecewise monotone interval
map $f$ a semiconjugated map $\tilde f$ which posesses a
countable Markov partition into intervals (the so-called
Hofbauer tower). Using this tool, Hofbauer and Keller [HK]
proved in 1984 that Theorem~2 from Section~1  holds even if
the initial interval map $f$ does not not admit a finite Markov
partitions into intervals where it is monotone.

The next generalisation consists in allowing $1/|f'|$ (or a more
general weight $g$) to be of bounded variation (sometimes, an addditional
assumption of continuity is used). The idea here is of course to
approach $g$ by a locally constant weight, where locally constant
changes meaning as the initial partition is refined by the dynamics
by considering finite intersections $\cap_{i=0}^M f^{-k} I_{t_k}$.
(Here, instead of assuming that $\widehat R < 1$ it is enough to
suppose that the initial partition is  generating, i.e. that
the diameter of
$\cap_{i=0}^\infty f^{-k} I_{t_k}$ is zero.) The corresponding
version of Theorem~2 was published in 1990 by Baladi and Keller,
a different proof (based on a slightly different ``Markov''
philosophy) is contained in the book of Ruelle [Ru2].

Piecewise injective (and piecewise expanding) maps have been studied
also in higher dimensions. The spectral theory of the transfer operator
is more technical (for example it is not obvious which Banach space to use!).
The survey [Ba3] contains references to the results of Saussol, Buzzi,
Tsujii, and others. A version of the Hofbauer tower can be constructed,
and Buzzi and Keller [BK] recently used it to prove an analogue
of Theorem~2 from Section~1 in the case when $f$ is piecewise
affine, piecewise expanding (in higher dimensions), and not
necessarily Markov.

The Hofbauer tower, or variants of it, has also been successful
to prove a version of Theorem~2 for one-dimensional (quadratic  e.g.)
maps with critical points (Keller-Nowicki) i.e. $c$ so that
$f'(c)=0$.

It must be noted that the approach we just described is a little
heavy to implement.

\smallskip

{\bf The Markov approach for smooth ``hyperbolic'' maps}

A different class of problems is given by $C^r$
($r > 1$) maps on compact
manifolds which are assumed to be uniformly expanding
or uniformly hyperbolic. There, classical results of Bowen, Sinai, Ruelle,
and others,
guarantee the existence of a finite Markov partition. (The definition
of Markov is slightly more
involved in the hyperbolic case, it also guarantees semiconjugacy with
a SFT.)
So the ``only'' problem here is to approach nonconstant weights $g$
(such as $|\det Df|^{-1}$) by locally constant weights.
For Holder $g$, this was done by Ruelle, Pollicott, and Haydn who
proved a version of Theorem~2 from Section~ 1 for hyperbolic diffeomorphisms.
In the case of higher smoothness $r \ge 2$ and $g$ at least $C^1$),
the natural object is in fact not the weighted zeta function, but
a weighted dynamical determinant of the type
$$
d(z)=\exp-\sum_{n=1}^\infty {z^n \over n}
\sum_{x \in \Fix f^n} {\prod_{k=0}^{n-1} g(f^k(x))\over
|\det (Df^{n}(x) -\Id)|}\, .
$$
It is possible to express the dynamical zeta function as an alternated product
of dynamical determinants. If $f$ and $g$ are $C^\omega$ (real-analytic)
then results from Cauchy and Grothendieck can be applied, and Ruelle [Ru0] showed
in 1976 that $d(z)$ is entire, and $\zeta_{f,g}(z)$ is meromorphic if $f$
is additionally assumed to be (locally) uniformly expanding. The hyperbolic case
produces serious additional difficulties, and the fact that
$\zeta(z)$ is meromorphic in the whole complex plane
was proved much more recently by Rugh [Rug] in dimension ~2 and
Fried [Fr] in general. The case of $C^r$ (non analytic) data was treated
by Ruelle [Ru1] for (locally) expanding maps. The hyperbolic $C^r$ case
is still partially mysterious, despite an important breakthrough by Kitaev [Kit].

Similar results exist for continuous-time dynamics (flows and semi-flows).

\medskip

{\bf The original Milnor and Thurston formula}

Let us now return to the situation of a continuous piecewise monotone
transformation
$f$ of a compact interval. Milnor and Thurston do not make any additional
assumptions
in [MT], so that the sets $\Fix f^n$ may have infinite
(even uncountable) cardinality.
However, if $f^n$ is decreasing on an interval $J$, it may have at most
one fixed point in $J$. Each set
$$
\Fix^- f^n =\{ x \in I \mid f^n (x)=x \, \,  ,
f^n \hbox{ is decreasing in a neighbourhood of } x \}
$$
therefore has finite cardinality,
and it makes sense to define a ``negative zeta
function''
$$
\zeta^-(z)=\exp \sum_{n=1}^\infty 2 \cdot \# \Fix^-(f^n)
$$
(the naive idea is that doubling ``negative''  fixed points makes up
for the ``forgotten'' ``positive''  fixed points -- a more precise
interpretation, making use of Lefschetz signs, is presented in
the 1996 paper of Ruelle quoted in [BaR]).

Milnor and Thurston then introduce the {\it kneading matrix.} If
$f$ has $N$ maximal intervals of monotonicity, it is an $(N-1)\times (N-1)$
matrix with coefficients power series in $z$,
the coefficients of which belong to $\{-1,0,+1\}$.
The $i$th line of this matrix is
$$
 {\theta(a_i+)(z) - \theta(a_i-)(z) \over 2} , \, i=1, \ldots, n-1\, ,
$$
where (we use the notation $\epsilon_f$ from Section ~1) the
kneading coordinate $\theta(x)(z)$ is the power series
$$
\theta(x)(z)=\sum_{k=0}^\infty z^k \prod_{j=0}^{k-1}\epsilon_f(f^j (x))
\cdot \alpha (f^k(x)) \, ,
$$
where $\alpha(y)$ is the $N-1$-tuple
$$
(\sgn(y-a_j), j=0, \ldots, N) \, ,
$$
with
$$
\sgn(\xi)=\cases
-1 & \xi < 0 \, , \cr
0&\xi=0\, ,\cr
1&\xi > 0 \, .
\endcases
$$

One of the key results in [MT] is the following remarkable equality:
\proclaim{Theorem 0 (Milnor-Thurston identity)}
$$
\zeta^-(z)={1 -z(\epsilon(a_0+)+\epsilon(a_N-)/2)\over \det(1+D(z))}\, .
$$
\endproclaim

As an immediate consequence, the negative zeta function is meromorphic
in the unit disc.

The above result is extremely beautiful, but (for the moment) a bit
mysterious. The proof of Milnor and Thurston (a homotopy argument
involving
the bifurcations of a path $f_t$ of piecewise monotone maps between $f=f_1$
and a ``trivial'' map $f_0$ having the same intervals of monotonicity as
$f$ and whose graph is strictly under the diagonal) does not give any
insight on ``why'' Theorem~1 holds. Note also that it is not clear
how to introduce weights in the negative zeta functions, and that
there is no spectral interpretation of the zeroes or poles of
the dynamical zeta function. The purpose of the remainder of Section~2
is to describe the one-dimensional kneading theory which addresses these
points.

\remark{Exercise 0}
Check that ${\theta(a_i+)(0) - \theta(a_i-)(0) \over 2}$
is indeed the vector $(0, \ldots, 0, 1, 0, \ldots , 0)$
where the $1$ is at the $i$th position.  Prove
the Milnor-Thurston identity for $f_0$.
\endremark

\smallskip
\subhead 2.2 The setting -- Essential spectral radius
\endsubhead

It will turn out to be more convenient to consider a slightly
more general setting, allowing the $\psi_i$ to be ``independent''
local homeomorphisms (in particular, with the possibility that
$\Imm \psi_i \cap \Imm \psi_j \ne \emptyset$ for $i\ne j$)
instead of the local inverse branches of a piecewise monotone
interval map $f$ as was the case until now.

\noindent {\bf The data}
\smallskip

We fix a compact interval $I \subset \real$ and a finite
set $\Omega$ of indices. (The restriction to finite
$\Omega$ is mostly for convenience and countable or
even uncountable index-sets endowed with a positive,
not necessarily finite, measure can also be used.
See [Ru3] and the exercises and remarks below.)

For each $\omega\in \Omega$, we take a nonempty open
subinterval $I_\omega$ of $I$ and a (local) homeomorphism
$$
\psi_\omega : I_\omega \to \psi_\omega (I_\omega) \, ,
$$
assuming $\psi_\omega (I_\omega) \subset I$ and
setting $\epsilon_\omega=+1$ if $\psi_\omega$ preserves
orientation and $\epsilon_\omega=-1$ otherwise.
We also consider a function $g_\omega : \real \to \complex$
satisfying:
\roster
\item
$g_\omega$ is supported in $I_\omega$,
\item
$\var g_\omega < \infty$,
\item
$g_\omega$ is continuous (on $I$, say).
\endroster

The third assumption allows us to use the Leibniz formula.
Taking the $\psi_\omega=\psi_j$ to be the local
inverse branches of a piecewise monotone $f$
(with $I_j=f(a_j, a_{j+1})$), and
$g_j = \chi_{I_j} |f' \circ \psi_j|^{-1}$ this assumption is not satisfied by all examples in Section~1,
but it can be essentially weakened, see [Go].

\remark{The dual system (Exercise 1)}
Given $\Omega$, $\psi_\omega$ and $g_\omega$
as above we may introduce a new ``dual'' system by setting
$\widehat I_\omega=\psi_\omega (I_\omega)$, $\hat \psi_\omega=
\psi_\omega^{-1}$, and $\hat g_\omega = \chi_{\widehat I_\omega} \epsilon_\omega
g_\omega \circ \psi_\omega^{-1}$. Check that this new system
satisfies all of the conditions (in particular
$\hat g_\omega$ is continuous contrary to what its
expression might suggest). Abusing notation, we
shall write $\hat g_\omega = \epsilon_\omega
g_\omega \circ \psi_\omega^{-1}$.
\endremark

\smallskip
\noindent {\bf The transfer operators}
\smallskip

\definition{Definition of the transfer operators}
We associate to our data two transfer operators acting either on the Banach
space $L^\infty$ of bounded functions (modulo functions which vanish except
on an at most countable set), or on the Banach space $\BB=BV(I)/\NN$:
$$
\eqalign
{
\MM \varphi &= \sum_\omega g_\omega \cdot \varphi \circ \psi_\omega \, ,\cr
\widehat \MM \varphi &
=\sum_\omega \hat g_\omega \cdot \varphi \circ \hat \psi_\omega=
\sum_\omega
\epsilon_\omega g_\omega\circ \psi_\omega^{-1} \cdot
\varphi \circ \psi_\omega^{-1}\, .}
$$
We shall use the notation
$$
\widehat \RR = \limsup_{n \to \infty}
\sup_{\varphi\, , \sup|\varphi| \le 1} \sup_{I} | \widehat \MM^n \varphi|^{1/n} \, .
$$
\enddefinition

\remark{Exercise 2}
Show that both $\MM$ and $\widehat \MM$ are bounded on
both Banach spaces considered. Show that the $\limsup$ defining $\widehat \RR$ is in
fact a limit and is the spectral
radius of $\widehat \MM$
on $L^\infty$. If $\Omega$ is countable,
find sufficient assumptions on the $\psi_\omega$ and
$g_\omega$ which imply that both operators are bounded on
both Banach spaces.
\endremark

\remark{($\widehat \MM$ as the dual of $\MM$)}
A priori, $\widehat \MM$ depends on the data $I$, $\Omega$, $\psi_\omega$,
$g_\omega$.  (Using partitions of unity, it is  easy to obtain different
data giving rise to the same operator.)
It is possible to show [BaRu]  that in fact it only depends on $\MM$
as an operator (on $\BB$, say) and not on the representation
of $\MM$ given by the $g_\omega$ and $\psi_\omega$.
For the moment we shall not need this
fact (we  just have to be aware that a preferred representation
must always be given, at least implicitly) but we slightly abuse terminology by
viewing $\widehat \MM$
as a dual of $\MM$ (see also the following exercise).
\endremark

\remark{Exercise 3}
Check (using each time the ``obvious'' representation
of the transfer operator) that $\widehat{\widehat \MM}=\MM$
and that $\widehat {\MM_1 \MM_2}=\widehat{\MM_2}\widehat {\MM_1}$
for all transfer operators $\MM$, $\MM_1$, $\MM_2$.
\endremark

\remark{Exercise 4}
If the $\psi_\omega$ are the local inverse
branches of a piecewise monotone interval
map and  $g_\omega = \chi_{I_\omega} g\circ \psi_\omega$ for
a single function $g$, find a simpler expression for
$\widehat \MM$ (and check that it is compatible with the notation
used in Section~1).
\endremark

For convenience, let us now introduce terminology that we have avoided
until now (see [Ba2] for more):

\definition{Definition of essential spectral radius}
Let $L:B\to B$ be a bounded linear operator on a Banach space $B$.
The {\it essential spectral radius} $\rho_\ess(L)$ of $L$ is
$$
\eqalign{
\rho_\ess(L)=
\inf \{ \rho > 0 \mid&
\hbox{ if } \lambda \in \sp(L) \hbox{ and }\cr
&\, \, |\lambda| > \rho \Longrightarrow \lambda
\hbox{ is an isolated eigenvalue of finite multiplicity} \}\, .
}
$$
\enddefinition

In other words, outside of the disc of radius $\rho_\ess$ the
spectrum is just like the spectrum of a compact
operator. If one can prove that the essential spectral
radius is strictly smaller than the spectral radius,
one often says that the operator is quasicompact.

\proclaim{Theorem 1 (Bound on $\rho_\ess(\MM)$)}
Let $\MM$ be as above. Then the essential spectral radius of
$\MM$ acting $\BB=BV/\NN$ is at most equal to $\widehat R$.
\endproclaim

\remark{Remark}
It is possible to show that the spectral radius of $\MM$
on $\BB$ is at most $\max(R,\widehat R)$ where $R$
is the spectral radius of $\MM$ acting on $L^\infty$.
The situation discussed at the end of Section~1 involved
the quasicompact case where $\widehat R < R$.
One can easily construct examples of data so that
$R < \widehat R$ (just take $\Omega$ a singleton and
$\psi_\omega$ a linear expansion) or $\widehat R=R$
(in the case where $\psi_\omega$ is the identity, e.g.).
Note that it is possible to obtain {\it lower bounds} on the
essential spectral radius, but this is much more tricky,
since we have to control sums of nonnecessarily positive
numbers.
(See [Go].)
\endremark

\demo{Proof of Theorem 1}
This can be proved just like Theorem~1 of Section ~1.
Checking it is a good exercise.\qed
\enddemo

\subhead 2.3 Sharp traces and sharp determinants
\endsubhead

We next define a (formal) trace and its associated
(formal) determinant  for transfer operators,
the ``sharp trace'' $\Tr^\# \MM$ and ``sharp determinant''
$\Det^\#(1-z\MM)$. We shall immediately prove some
of their basic properties, but it is only in the next section that
we shall introduce the so-called ``kneading operators'' which will allow
us to show (in \S~ 2.5)
that the zeroes of $\Det^\#(1-z\MM)$ in a suitable disc
describe some of the inverse eigenvalues of $\MM$ on $\BB$.

It is possible to show [BaRu] that the following definition only
depends on $\MM$ as an operator on $\BB$ (instead, we always
assume that a preferred representation is given):
\definition{Definition (Sharp trace and sharp determinant)}
Let $\MM$ be a transfer operator associated to data $\psi_\omega$,
$g_\omega$ as in \S ~2.2. Then we write:
$$
\Tr^\# \MM = \sum_\omega \int {1 \over 2} \sgn(\psi_\omega(x)-x) \,
dg_\omega(x)\, ,
$$
and (as a formal power series in $z$)
$$
\Det^\#(1-z\MM)=
\exp - \sum_{n=1}^\infty {z^n\over n}
\Tr^\# \MM^n \, .
$$
\enddefinition

The expression for the trace is well-defined since each $dg_\omega$
is a complex measure and the integrand is a bounded function.
Note for future use that $\sgn(\psi_\omega x -x)$ is a
function of bounded variation.

\remark{Exercise 5}
Show that $d\sgn$ is twice the Dirac mass at $0$. (See also
the proof of Lemma~2 in \S ~1). Rewrite the expression for $\Tr^\# \MM$
if each $g_\omega$ is $C^1$.
\endremark

\smallskip

\noindent Let us first prove an easy but very useful result:
\proclaim{Lemma 1 (``Functional relation'')}
$\Tr^\#\widehat \MM = - \Tr^\# \MM$. In particular
$$
\Det^\#(1-z\widehat \MM) = {1 \over \Det^\#(1-z\MM) } \, .
$$
\endproclaim

\demo{Proof of Lemma 1}
The proof is based on the change of variable formula
(Lemma ~2 from Section ~1):
$$
\eqalign
{
\Tr^\# \widehat \MM &= \sum_\omega \int \epsilon_\omega {1\over 2}
\sgn(\psi_\omega^{-1}(x)-x) \, d(g_\omega \circ \psi_\omega^{-1}(x))\cr
&=\sum_\omega \int {1\over 2}
\sgn(y-\psi_\omega(y)) \, d(g_\omega (y))\cr
&=-\Tr^\# \MM \, .
}
$$
The second claim is left as an exercise.
\qed
\enddemo

The next lemma requires more effort, but it gives legitimacy
to the ``trace'' terminology (not yet in the sense that the
trace is related to the eigenvalues, however):

\proclaim{Lemma 2 (Trace property)}
$$
\Tr^\# (\MM_1 \MM_2) = \Tr^\# (\MM_2 \MM_1) \, .
$$
\endproclaim

\remark{Exercise 6}
Show that Lemma~2 implies that we may perform the usual algebraic
manipulations on the sharp determinants i.e.:
$$\eqalign{
&\Det^\#(1-z\MM_1 \MM_2)= \Det^\#(1-z\MM_2 \MM_1)\cr
&\Det^\#(1-z\MM_1) \Det^\#(1-z\MM_2)=
\Det^\#(1-z(\MM_1+\MM_2 + z \MM_1 \MM_2))\, .
}
$$
\endremark

\demo{Proof of Lemma 2} By linearity of the sharp trace, it
is enough to prove the lemma in the case where $\Omega_1$
and $\Omega_2$ are both singletons i.e.
$$
\MM_1 \varphi = g_1 \cdot \varphi \circ \psi_1 \, , \quad
\MM_2 \varphi = g_2 \, \varphi \circ \psi_2 \, .
$$
First assume that $\psi_2 \psi_1$ is increasing, and let
$\epsilon= \pm 1$ depending on whether $\psi_1$ and $\psi_2$
are increasing or decreasing. Since
$\psi_1$ and $\psi_2$ are continuous, the set
$\{ x \, : \, \psi_2 \psi_1 x \ne x \}$ is the union of
at most countably many open
intervals $(c_i, d_i)$. Correspondingly,
$\{ y \, : \, \psi_1 \psi_2 y \ne y \}$ is the union of intervals
$(c'_i, d'_i)$ where
$$
c'_i = \psi_1 c_i = \psi_2^{-1} c_i \, ,\quad
d'_i = \psi_1 d_i = \psi_2^{-1} d_i \, ,
$$
if $\epsilon = 1$ and
$$
c'_i = \psi_1 d_i = \psi_2^{-1} d_i \, ,\quad
d'_i = \psi_1 c_i = \psi_2^{-1} c_i \, ,
$$
if $\epsilon = -1$. If $\sigma_i$ is the sign of $\psi_2 \psi_1
x-x$ on $(c_i,d_i)$, then $\sigma'_i = \epsilon \sigma_i$ is the sign
of $\psi_1 \psi_2 y - y$ on $(c'_i, d'_i)$.

We have
$$
\eqalign
{
\Tr^\# \LL_1 \LL_2 & =
\int d(g_1(x)g_2(\psi_1(x))
{1\over 2} \sgn (\psi_2 \psi_1 x -x) \cr
&={1\over 2} \sum_i \int_{c_i}^{d_i} d(g_1(x) g_2(\psi_1(x)) \sigma_i \cr
&={1 \over 2 }
\sum_i \sigma_i
\bigl [ g_1 (d_i) g_2 (\psi_1 d_i) -
g_1 (c_i) g_2 (\psi_1 c_i) \bigr ]
\cr
&={1 \over 2}
\sum_i \sigma_i \epsilon
\bigl [
g_1 (\psi_2 d'_i) g_2 (d'_i) -
g_1 (\psi_2 c'_i) g_2 (c'_i) \bigr ] \cr
&=
{1 \over 2}
\sum_i
\sigma'_i
\bigl [
g_2(d'_i) g_1 (\psi_2 d'_i) -
g_2 (c'_i) g_1 (\psi_2 c'_i) \bigr ] \cr
&=
\int
d(g_2(y) g_1(\psi_2(y)) {1\over 2}
\sgn (\psi_1 \psi_2 y - y) \cr
&= \Tr^\# \LL_2 \LL_1 \, . \cr
}
$$

If $\psi_2 \psi_1$ is decreasing, either it has
no fixed point and $\psi_1 \psi_2$ has
no fixed point either, or it has a unique fixed point
$c$ and
$$
c'=\psi_1 c = \psi_2^{-1} c
$$
is the unique fixed point of $\psi_1 \psi_2$. Then
$$
\eqalign
{
\Tr^\# \LL_1 \LL_2 &=
g_1(c) g_2(\psi_1 c) \cr
&= g_2(c') g_1 (\psi_2 c') \cr
&= \Tr^\# \LL_2 \LL_1
}
$$
concluding the proof.\qed
\enddemo

\remark{Exercise 7 (Sharp determinant as a Lefschetz weighted zeta function)}
Let $\MM$ be given by data such that for each $n\ge 1$
and each $\omega_1, \ldots, \omega_n$ with
$$g_{\omega_n}\circ\psi_{\omega_{n-1}} \circ \ldots \circ \psi_{\omega_1}
\cdots g_{\omega_2} \circ \psi_{\omega_1} \cdot  g_{\omega_1}$$
well-defined and nonzero, the set of points fixed by
$\psi_{\omega_n} \circ \ldots \circ \psi_{\omega_1}$ is finite.
Assume furthermore that each $\psi_{\omega}$ is a local diffeo.
Show that
$$
\Tr^\#\MM=\sum_\omega \sum_{x: \psi_\omega(x)=x}
g_\omega(x) \sgn(1-\psi_\omega'(x))\, .
$$
Write the analogous expressions for $\Tr^\#\MM^n$
and $\Det^\#(1-z\MM)$.
\endremark

\subhead 2.4 Kneading operators and the key equality
\endsubhead

We consider data $\psi_\omega$ and $g_\omega$
as in \S~2.2. We finally introduce the weighted analogues of the Milnor and Thurston kneading
matrices and prove the corresponding ``weighted Milnor-Thurston''
identity between determinants.

Let us first consider a slightly {\it simpler case} where the notation
is less heavy: assume additionally that the $g_\omega$ are $C^1$
and the $\psi_\omega$ are local diffeomorphisms.
We then introduce two auxiliary transfer operators (we shall see later
where they act, we can take e.g. $\BB$, $L^\infty$ or $L^2(Leb)$):
$$
\eqalign
{
\NN \varphi&=\sum_\omega g'_\omega \cdot \varphi \circ \psi_\omega \, , \cr
\widehat \NN \varphi&=\sum_\omega \epsilon_\omega
g'_\omega\circ \psi_\omega^{-1} \cdot
(\psi_\omega^{-1})' \cdot \varphi \circ \psi_\omega^{-1} \, . \cr
}
$$
We also introduce a convolution operator (mapping e.g.  $L^\infty$ functions
with compact support to $L^\infty$ functions):
$$
\eqalign
{
\SS \varphi (x)&=\int {1\over 2} \sgn(x-y) \, \varphi(y)\, dy\cr
&=\int {1\over 2}  \sgn(y) \, \varphi(x-y)\, dy\, .\cr
}
$$
Finally we define the kneading operators $\DD(z)$ and $\widehat \DD(z)$ in the
sense of formal power series with operator coefficients
(using the shorthand notation $(\Id-z\MM)^{-1}$
for the formal power series $\sum_{k\ge 0} z^k \MM^k$):
$$
\DD(z)= z\NN (\Id-z\MM)^{-1} \SS \, , \,
\widehat \DD(z)=z\widehat \NN (\Id-z\widehat \MM)^{-1} \SS \, .\tag{1}
$$
The description of $\DD(z)$ and $\widehat \DD(z)$
below as kernel operators will show that the series above define
bounded operators on $L^2(d\mu)$ whenever $1/z \notin \sp \MM$
(acting on $\BB$) respectively $1/z \notin \sp \widehat \MM$
(acting on $L^\infty$).

\smallskip

{\bf Although $\SS$ can be viewed as a  bounded operator from
$L^2(Leb,I)$ to itself,
it  does not map $L^2(Leb)$ into itself
boundedly, and $(\Id-z\MM)^{-1}$ is not necessarily bounded on $L^2(Leb)$ or $L^2(Leb,I)$
even if $1/z \notin \sp(\MM|\BB)$.}

\remark{Exercise 8}
The notation $\SS$ is not strictly compatible with the notation
from Section~1. Show that the operator $\SS$ from Lemma~0 in
Section ~1, when acting on (the densities of) Radon measures which are absolutely
continuous with respect to Lebesgue, can be written as
$$
\SS \varphi(x)=\int {1\over 2} (\sgn (y)+\sgn(x-y)) \, \varphi(y) \, dy \, .
$$
\endremark

\smallskip
In the {\it general case} (i.e. $g_\omega$ is continuous and of bounded variation
while $\psi_\omega$ is a local homeomorphism), it is useful to associate
a finite nonnegative measure  to our data:
$$
\mu = \sum_\omega |dg_\omega| + \sum_\omega |d (g_\omega \circ \psi_\omega^{-1})| \, .
$$
This measure is constructed in such a way as to guarantee that
the Radon-Nikodym derivatives $dg_\omega /d\mu$ and
$d (g_\omega \circ \psi_\omega^{-1})/d\mu$  exist and are bounded.
We can now redefine the auxiliary transfer operators:
$$
\eqalign
{
\NN \varphi&=\sum_\omega {dg_\omega \over d \mu}\cdot \varphi \circ \psi_\omega \, , \cr
\widehat \NN \varphi&=\sum_\omega \epsilon_\omega
{d(g_\omega\circ \psi_\omega^{-1}) \over d\mu}\cdot \varphi \circ \psi_\omega^{-1} \, . \cr
}
$$
Similarly, we redefine the convolution operator as:
$$
\SS \varphi (x)=\int {1\over 2} \sgn(x-y) \, \varphi(y)\, d\mu\, .
$$
Finally, $\DD(z)$ and $\widehat \DD(z)$ are defined as in
equation \thetag{1} above.
Note that $\SS$ is bounded from $L^2(d\mu)$ to $L^2(d\mu)$
(restricting $x$ to the compact support of $\mu$ in the left-hand-side),
but this will not be very useful since the resolvent $(\Id-z\MM)^{-1}$
is not bounded on $L^2(d\mu)$ for all $1/z \notin \sp(\MM)$.
What we will do next is notice that $\DD(z)$ and
$\widehat \DD(z)$ are kernel operators and examine their
kernels for appropriate $z$, i.e.  write
$$
\DD(z) \varphi(x)= \int_I K^z(x,y) \, \varphi(y) \, d\mu(y) \, .
$$
This is straightforward: we have, formally
$$
\eqalign
{
\DD(z) \varphi (x)&= z \NN (\Id-z\MM)^{-1} \SS \varphi(x)\cr
&= \int [z \NN_x (\Id-z\MM)^{-1}_x {1\over 2} \sgn (\cdot -y)] (x)\,
\varphi(y)\, d\mu(y)\cr
&= \int \sum_\omega z {d g_\omega  \over d\mu}(x) (\Id-z\MM)^{-1}_x
{1\over 2} \sgn (\cdot -y)] (\psi_\omega x)\,
\varphi(y)\, d\mu(y)\, .\cr
}
$$
(The index in $\NN_x$ or $(\Id -z\MM)^{-1}_x$ is here to emphasize on which
variable the transfer operator is acting.)
Since $\sgn (\cdot -y)/2$ is of bounded variation (uniformly) for each $y$,
it is clear from the above expression that if $1/z \notin \sp(\MM)$
(on $\BB$) then the kernel  $K^z(\cdot,y)$ is bounded 
uniformly in $y$, therefore a bounded function of $x \in I$ and $y$ in
the support of $\mu$ and since $L^\infty(I\times I)\subset 
L^2(d\mu\times d\mu)$
the kernel is in $L^2(d\mu\times d\mu)$.
Similarly
$$
\eqalign
{
\widehat \DD(z) \varphi (x)&= z \widehat \NN (\Id-z\widehat \MM)^{-1} \SS \varphi(x)\cr
&= \int \sum_\omega z {d \widehat g_\omega  \over d\mu}(x) (\Id-z\widehat \MM)^{-1}_x
{1\over 2} \sgn (\cdot -y)] (\widehat \psi_\omega x)\,
\varphi(y)\, d\mu(y)\, ,\cr
}
$$
so that the kernel of $\widehat \DD(z)$ is in $L^2(d\mu\times d\mu)$
for all $1/z \notin \sp\widehat \MM$, where we  consider
$\widehat \MM$ acting on bounded functions.
It follows that

\proclaim{Lemma 3}
For all $1/z \notin \sp(\MM)$ (on $\BB$) the operator
$\DD(z)$ is a compact operator when acting on $L^2(d\mu)$, it
is in fact Hilbert-Schmidt.
For all $1/z \notin \sp(\widehat \MM)$ (on $L^\infty$) the operator
$\widehat \DD(z)$ is a compact, in fact Hilbert-Schmidt operator when acting on $L^2(d\mu)$.
\endproclaim

In particular, the spectra of both kneading operators for $|z|<1/\widehat R$ 
and
$1/z\notin \sp(\MM|\BB)$ consist in eigenvalues of finite 
multiplicity which can only
accumulate at $0$, and which are the zeroes of  entire
functions, the regularised determinants (of order two)
$$
\lambda \mapsto \Det_2(1-\lambda \DD(z)) \, ,
\quad \lambda \mapsto \Det_2(1-\lambda \widehat \DD(z)) \, .
$$
We refer to Appendix A for the part of the theory of Hilbert-Schmidt operators
that we  use.
\smallskip

Before finally stating the weighted equivalent
of the Milnor and Thurston formula, let
us propose an exercise which  makes the link between the
notation used here and the original definitions in [BaRu].

\remark{Exercise 9}
If $\mu$ is a finite nonnegative measure on $I$ and
$K\in L^2(d\mu\times d\mu)$, we associate to an operator
$$
D \varphi(x)=\int K(x,y)\, \varphi(y)\, d\mu(y) \, ,
$$
acting on $L^2(d\mu)$ another operator on
$L^2(d\mu)$, noted $D^*$, by setting
$$
D^* \varphi (x)=\int K(y,x) \, \varphi(y) \, d\mu(y) \, .
$$
\roster
\item
Show that for every $m\ge 2$, the kernel $K_m(x,y)$
of $D^m$ coincides with the kernel $K^*_m(x,y)$
of $(D^*)^m$ on the
diagonal $x=y$. (If one has a preferred representant for
$K(x,y)$, then the statement also makes sense for $m=1$,
otherwise the fact that the diagonal has zero measure in
the square $I\times I$ causes problem: understanding this is
necessary to do the exercise! See also Lemma~6 in Section 2.)
\item
Check that the kernel of $\DD(z)^*$ is
$$
\sum_\omega z {dg_\omega \over d \mu} (y)
[(\Id -z \MM)^{-1}_y {1\over 2}\sgn(\cdot - x)] \psi_\omega(y)\, .
$$
\endroster
\endremark

\smallskip

{\bf The weighted Milnor and Thurston identity}

Since $\DD(z)$ and $\widehat \DD(z)$ are Hilbert-Schmidt on $L^2(d\mu)$ for each
$z< 1/\widehat R$ and $1/z \notin \sp(\MM)$ (on $\BB$), the theory
in Appendix~A allows us to introduce formal determinants:
$$
\eqalign
{\Det_*(\Id+\DD(z))&=
\exp (+\int_I K^z(x,x) \, d\mu(x)) \cdot \Det_2(\Id+\DD(z))\cr
\Det_*(\Id+\widehat \DD(z))&=
\exp (+\int_I \widehat K^z(x,x) \, d\mu(x)) \cdot \Det_2(\Id+\widehat \DD(z))\, ,\cr
}
$$
where we also used the fact that the kernels $K^z(x,y)$
and $\widehat K^z(x,y)$ of $\DD(z)$ and
$\widehat \DD(z)$ are bounded functions well-defined almost
everywhere in the diagonal $x=y$.

The first important result of this section is:
\proclaim{Theorem 4 (Weighted Milnor-Thurston identity)}
In the sense of formal power series, we have:
$$
\Det^\#(1-z\MM) = \Det_*(1+\widehat \DD(z))=
{1 \over \Det_*(1+\DD(z)) }\, .
$$
\endproclaim

To prove Theorem ~4, we shall express $\Det_*(\Id+\DD(z))$
as an exponential of a sum of ``Fredholm-type'' traces, i.e.
averages of the kernels of  $\DD^n(z)$ on the diagonal
of $I\times I$. For this, the following  lemma is essential:

\proclaim{Lemma 5 (Trace and Fredholm trace for $L^2$ kernels)}
Let $\mu$ be a finite nonnegative Borel measure on $I$.
Let $A : L^2(d\mu)\to L^2(d\mu)$ be a Hilbert-Schmidt
operator  described by (an $L^2 (d\mu \times d\mu)$
kernel $K :I \times I\to \complex$:
$$
A\varphi(x)= \int K(x,y) \varphi(y) \, d\mu(y)\, .
$$
(In particular, $A^2\in \SS^1$.)
Then $\Tr A^2 = \int K(x,y) K(y,x) \, d\mu (y) \, .$
\endproclaim

\demo{Proof of Lemma 5}
Combine the Lidskii theorem in the Appendix with exercise~49 of 
Chapter XI in [DS2].
\qed
\enddemo

\remark{Consequence of Lemma 5}

If $|z|< 1/\widehat R$ and $1/z\notin \sp(\MM)$ (on $\BB$)
then, writing $K^{z,n}(x,y)$ for the kernel of $\DD(z)^n$
acting on $L^2(d\mu)$, we have
$$
\Det_*(1+\DD(z)) =\exp-\sum_{n=1}^\infty {(-1)^n\over n}
\int_I K^{z,n}(x,x)\, d\mu(x) \, .
$$
In the proof of Theorem~4, we shall view
the above expression only in the sense of formal power series, and it
will be convenient to use the decomposition
$\DD(z)=z\NN \sum_{k=0}^\infty z^k \MM^k \SS$ and the
corresponding  expressions for the $K^{z,n}(x,y)$
as formal power series.
\endremark

Our final ingredient for the proof of Theorem~4 is
the following purely algebraic exercise on formal
traces and formal determinants:
\remark{Exercise 10 (Properties of formal determinants)}
Let $\AA$ be a vector space over $\complex$ which is a
subset of an algebra. We write
$\AA^\infty$ for the set $\{ \KK \in \AA \mid \KK^n \in \AA\, ,
\forall n \ge 1 \}$.
To any  function (the {\it formal trace})
$$
\widetilde {\hbox{tr}}\,  :\AA \to \complex \, ,
$$
we associate a {\it formal determinant}
$$
\widetilde {\hbox{det}} (\Id + \lambda \cdot) : \AA^\infty
\to \complex [[ \lambda]] \, ,
$$
by setting
$$
\widetilde {\hbox{det}}\, (\Id+\lambda \KK)
= \exp -\sum_{n=1}^\infty {(-\lambda)^n\over n} \widetilde {\hbox{tr}}\, \KK^n \, .
$$
Show that:
\roster
\item
If $\KK_1, \KK_2 \in \AA$ are such that $\KK_1 \KK_2 \in \AA^\infty$
and $\KK_2 \KK_1 \in \AA^\infty$, and
$$\widetilde {\hbox{tr}}\,( (\KK_1 \KK_2)^n )=
\widetilde {\hbox{tr}}\,( (\KK_2 \KK_1)^n )
$$
for all integer $n\ge 1$, then
$$
\widetilde {\hbox{det}}\, (\Id+\lambda \KK_1 \KK_2)=
\widetilde {\hbox{det}}\, (\Id+\lambda \KK_2 \KK_1)\, .
$$
\item
If, additionally, $\KK_1$ and $\KK_2$ are such that,
for each integer $m\ge 1$ (say, even) and every $m$-tuple $j_1, \ldots, k_m$
of nonnegative integers, we have
$$
\KK^m_{\vec j} = \KK_1^{j_1} \KK_2^{j_2} \cdots \KK_1^{j_{m-1}}\KK_2^{j_m} \in \AA
$$
and  $\widetilde {\hbox{tr}}\,  \KK^m_{\vec j}$ is well-defined
and invariant
under circular permutations of the $j_\ell$, then
$$
\widetilde {\hbox{det}}\, (\Id+\lambda \KK_1) \cdot
\widetilde {\hbox{det}}\, (\Id+\lambda \KK_2)=
\widetilde {\hbox{det}}\, ((\Id+\lambda \KK_1) (\Id+\lambda \KK_1))\, .
$$
\endroster
(Note that both conditions hold if the first one is true for $m=1$
and all $\KK_1$ $\KK_2$, if
$\AA$ is an algebra. In our application we will not have an algebra
strictly speaking.)
\endremark

\demo{Proof of Theorem~4}
We shall use the notation
$$
\sigma(t) ={\sgn(t)\over 2} \, .
$$
In order to apply Exercise~10, we must formalise the operator spaces
which appear and show that the commuting property holds.
For this, we set $\AA^S$ to be the vector space generated
by  finite products
of operators $\MM$, $\NN$, and $\SS$ such that there is
at least one factor $\SS$, at least an $\MM$ or an $\NN$
between any two  $\SS$ factors,
and at least one factor $\SS$
between any two $\NN$ factors.  By definition $\DD(z) \in \AA^S[[z]]$.

The following two properties will allow
us to invoke Exercise~10:
\smallskip
\noindent{\bf(I) Properties of the kernel}

If $\KK \in \AA^S$ then there exists $K\in L^2(\mu\times \mu)$
so that (on $L^2(d\mu)$, say)
$$
\KK \varphi(x)=
\int K(x,y) \varphi(y) \, d\mu (y) \, ,
$$
and, additionally, $K(x,y)$ is a (finite) linear combination
of expressions
$$
h(x) \cdot \tilde h(y) \cdot \sigma(\psi(x)-\tilde \psi(y)) \, ,
$$
where
\roster
\item
$\psi$ and $\tilde \psi$ are homeomorphisms, or local homeomorphism
the support of which contain the supports of $h$, respectively $\tilde h$;
\item
$h$ is a linear combination of continuous
functions of bounded variations, which may be multiplied by a factor
$(dg_\omega/d\mu) \circ \psi'$ with $\psi'$
a (local) homeomorphism (with good support)
if there is a factor $\NN$ which
is not followed by a (post)\-com\-posi\-tion with an $\SS$;
its support is compact if the leftmost factor of
$\KK$ is not $\SS$;
\item
$\tilde h$ is bounded; its support is compact if the rightmost factor of
$\KK$ is not $\SS$;
\endroster

The kernel $K$ is not uniquely defined, but our proof of its existence
is constructive so that there is no ambiguity once the data
$g_\omega$, $\psi_\omega$ is given.
Note that property (I) says in particular that the trace
$$
\Tr^* \KK = \int K(x,x) \, d\mu(x) \, ,
$$
is well defined for each element of $\AA^S$. Renaming $\Tr^*$
(for the sake of uniform notation)
$\Tr^\#$ on this vector space, and setting $\XX$ to be the
algebra generated by (powers of) $\MM$, we may extend $\Tr^\#$
by linearity to $\AA=\XX [[z]]\oplus \AA^S[[z]]$.
All expressions appearing in the proof of Theorem~4 will
be in $\AA$.

\smallskip
\noindent{\bf(II) Commutations}
\roster
\item
If $\KK \in \AA^S$ then
$\Tr^\# (\MM \KK) = \Tr^\# (\KK \MM)$ and
$\Tr^\# (\NN \KK) = \Tr^\# (\KK \NN)$.
\item
If $\KK \in \AA^S$ and neither the leftmost nor the rightmost
factor of $\KK$ is $\SS$ then $\Tr^\# (\SS \KK)=\Tr^\# (\KK \SS)$.
\item
$\Tr^\# (\MM \NN) =\Tr^\# (\NN \MM)$
(in fact, we do not use this in the proof; recall also Lemma~2).
\endroster

\smallskip
\noindent The third property will be crucial in the proof:

\noindent{\bf(III) Naturality}

For each $m\ge 1$
$$
\Tr^\# ((\MM-\SS \NN)^m) = 0 \, .
$$
\smallskip

We next show how Theorem~4 follows from Exercise~ 10 and (I-II-III).
It suffices to show that $\Det^\#(\Id+\DD(z)) \Det^\#(\Id-z\MM)\equiv 1$:
$$
\eqalign
{
\Det^\#(\Id+\DD(z)) \Det^\#(\Id-z\MM)
&= \Det^\#(\Id+z \NN(\Id-z\MM)^{-1} \SS) \Det^\#(\Id-z\MM)\cr
&=\Det^\#(\Id+z \SS \NN(\Id-z\MM)^{-1} ) \Det^\#(\Id-z\MM)\cr
&=\Det^\#(\Id-z\MM +z \SS \NN ) =1 \, , \cr
}
$$
where we used the definition of $\DD(z)$ in the first equality,
(I) and (II \therosteritem {2}) (together with Exercise~10) in the second
one, and (I) and (II \therosteritem {1, 2})
(with Exercise ~10 again) in the third, and
(III) in the last.

\smallskip
It remains to check (I,II, III).

We shall prove (I) by induction on the number of factors, multiplying to
the left. The claim is obvious for $\KK= \SS$. For $\SS\MM$
and $\MM\SS$, we  compute:
$$
\eqalign{
\SS \MM\varphi(x)&=
\int \sigma(x-y) \sum_\omega g_\omega(y) \varphi(\psi_\omega(y)) \, d\mu(y)\cr
&=\int \sum_\omega \epsilon_\omega \chi_{\psi_\omega(I_\omega)}
\sigma(x-\psi_\omega^{-1}(z)) g_\omega(\psi_\omega^{-1}(z)) \varphi(z) \, d\mu(z) \, ;\cr
}
$$
and
$$
\MM \SS\varphi(x)=
\int \sum_\omega g_\omega(x) \sigma(\psi_\omega(x)-y) \varphi(y) \, d\mu(y) \, .
$$
The above computations also show that $\SS \NN$ and $\NN \SS$
have kernels with the desired properties.

Next, assuming that the kernel of $\KK$ in $\AA$ has the desired
properties, we consider $\MM\KK$, $\NN \KK$ (if there is at
least an $\SS$ postcomposed with the last
factor $\NN$ in $\KK$),
$$
\eqalign{
\MM \KK \varphi(x)
&= \int \sum_\omega g_\omega(x) K(\psi_\omega(x), y) \varphi(y)
\, d\mu(y) \, ,\cr
\NN \KK \varphi(x)
&= \int \sum_\omega {dg_\omega\over d\mu }(x) K(\psi_\omega(x), y) \varphi(y)
\, d\mu(y) \, ,}
$$
for which it is obvious that the induction hypotheses suffice.

Finally, we consider $\SS \KK$ (if the leftmost
factor of $\KK$ is not $\SS$):
$$
\eqalign{
\SS \KK \varphi(x)&=
\int \int \sigma(x-y) K(y,z) \varphi(z) \, d\mu(z) \, d\mu(y)\cr
&=\int \int \sigma(x-y) K(y,z)  \, d\mu(y) \, \varphi(z)\, d\mu(z)\cr
&= \int \int \sigma(x-y) h(y) \tilde h(z) \sigma(\psi(y)-\tilde \psi(z))
\, d\mu(y) \, \varphi(z) \, d\mu(z) \, .
}
$$
By induction, the support of $h$ is compact. We must study
$$
\int \sigma(x-y) h(y)  \sigma(\psi(y)-u)  \, d\mu(y)   \, ,
$$
where we wrote $u=\tilde \psi(z)$ for simplicity. If $h$ is a linear
combination of continuous functions of bounded variation, we use the
easily proved fact
that $\SS_{1\to 0} (h d \mu)$ is a continuous function of bounded
variation, where $\SS_{1\to 0}$ is the isomorphism between Radon
measures and $\BB$ from Lemma~0 in Section~1 (the isomorphism was called $\SS$
here, but the notation $\SS$ in the present Section~2 represents the
convolution operator $\SS=\SS_{0\to 0}$). If $\KK$ contains an $\NN$ factor
not postcomposed by any $\SS$, then $h$ may contain
terms of the form $h'(y) (dg_\omega/d\mu)(\psi' y)$ with $h'$
continuous and of bounded variation, and
$$\SS_{1\to 0} (h'{dg_\omega\over d\mu} \circ \psi' d \mu)=
\SS_{1\to 0} (\epsilon_{\psi'}\chi_{\psi'} h'\circ (\psi')^{-1}dg_\omega)
$$
is again a continuous function of
bounded variation. Therefore, we may use the Leibniz formula
(Lemma~1 from Section~1) in both cases:
$$
\eqalign{
\int \sigma(x-y)  \sigma(\psi(y)-u)  h(y) \, d\mu(y)
&= \int \sigma(x-y)  \sigma(\psi(y)-u) d \SS_{1\to 0}( h d\mu) (y)\cr
&=   \sigma(\psi(x)-u) \SS_{1\to 0}(h d\mu) (x)\cr
&\qquad- \epsilon_\psi \chi_\psi
\sigma(x- \psi^{-1}(u))   \SS_{1\to 0}(h d\mu) (\psi^{-1}(u)) \, ,
}
$$
where we used $d\sigma=\delta_0$
again (the assumptions on the support of $h$ imply that there is
no boundary term). Inspecting the above expression, we see that we have
performed the inductive step successfully.

\smallskip
Let us prove II(1).  Using the expression obtained previously for $\MM\KK$,
and a change of variables, we get
$$
\eqalign{
\Tr^\# \MM \KK&=
\int \sum_\omega K (\psi_\omega(y), y) g_\omega(y) \, d\mu(y) \cr
&= \int \sum_\omega \epsilon_\omega \chi_{\psi_\omega(I_\omega)}
K(x, \psi_\omega^{-1}(x)) g_\omega(\psi_\omega^{-1}(x)) \, d\mu(x) \, .
}
$$
Similarly, we have:
$$
\eqalign{
\KK\MM \varphi (y) &=
\int  K (y, z) \MM\varphi(z) \, d\mu(z) \cr
&= \int  K(y,z) \sum_\omega  g_\omega(z)  \varphi(\psi_\omega (z)) \, d\mu(z) \cr
&=\int \sum_\omega \epsilon_\omega \chi_{\psi_\omega(I_\omega)}
K(y, \psi_\omega^{-1}(x)) g_\omega(\psi_\omega^{-1}(x))
\varphi(x) \, d\mu(x) \, ,
}
$$
which gives $\Tr^\# \MM\KK=\Tr^\# \KK \MM$.
Since we did not use integration by parts, the same
computation yields $\Tr^\# \NN\KK=\Tr^\# \KK \NN$.

To show II(2), we first use the expression
for $\SS \KK$ to see that
$$
\Tr^\# \SS \KK=
\int \int \sigma(x-y) K(y,x) \, d\mu(y) \, d\mu(x) \, .
$$
On the other hand
$$
\eqalign{
\KK \SS \varphi(x)&=
\int  K(x,y) \SS \varphi(y) \, d\mu(y) \cr
&= \int \int K(x,y) \sigma(y-z) \varphi(z) \, d\mu(z) \, d\mu(y) \, ,}
$$
so that
$$
\Tr^\# \KK \SS = \int \int K(x,y) \sigma(y-x)  \, d\mu(x) \, d\mu(y) \, ,
$$
proving the claim.

The proof of II(3) goes along the lines of the
proof of Lemma~2 (the Leibniz formula is all right since there
is only a single factor $\NN$).
\smallskip
Finally, we check (III) by induction on $m$. For $m=1$,
using our formula for $\Tr^\# \KK \NN$ in the case
$\KK=\SS$ (so that $K(x,y)=\sigma(x-y)$),
we find by a double change of variable
$$
\eqalign{
\Tr^\# \SS \NN &=
\int \sum_\omega \sigma(x-\psi_\omega^{-1}(x)) {dg_\omega\over d\mu}
(\psi_\omega^{-1}(x) ) d\mu(x)\cr
&=\int \sum_\omega  \sigma(x-\psi_\omega^{-1}(x))
dg_\omega (\psi_\omega^{-1}(x) )
{d\mu \over d\mu \circ \psi_\omega^{-1}}(x) \cr
&=\int \sum_\omega \sigma(\psi_\omega(y)-y) \, dg_\omega(y) \cr
&=\Tr^\# \MM \, .
}
$$
Next, we set
$$
\widetilde \MM= \MM - \SS \NN \, .
$$
We have just seen that $\Tr^\# \widetilde \MM=0$ and we want to
show that
$\Tr^\# (\widetilde \MM)^m=0$ for all $m\ge 1$. For
this, it suffices to show that
$$(\widetilde \MM)^m=\widetilde {\MM^m}$$
(where the notations are self-explanatory i.e.
$\widetilde{\MM^m}=\MM^m - \NN_m$ where $\NN_m$ is associated to
$\MM^m$ via an appropriate auxiliary measure $\mu_m$).
Indeed, the case $m=1$ applied to $\MM^m$ would give the claim.
To prove the above naturality statement, it is enough
(by density of $\BB$ in $L^2(d\mu)$) to show that
$(\widetilde \MM)^m\varphi=\widetilde {\MM^m}\varphi$
for each $\varphi$ of bounded variation.
Let us rewrite $\SS \NN$ on $\BB$,
integrating by parts:
$$
\eqalign{
\SS \NN \varphi(x)&=
\int \sum_\omega \sigma(x-y) \varphi(\psi_\omega(y)) \,  dg_\omega(y)\cr
&= \sum_\omega g_\omega(x) \varphi(\psi_\omega(x))\cr
&\qquad- \int \sum_\omega \sigma(x-y) g_\omega(y) d(\varphi \circ \psi_\omega)(y)\cr
&=\MM \varphi(x) -
\int  \sum_\omega  \epsilon_\omega \chi_{\psi_\omega(I_\omega)}
\sigma(x-\psi_\omega^{-1}(z))) g_\omega(\psi_\omega^{-1}(z)) d\varphi(z)\cr
&= \MM \varphi(x) -\NN_{1\to 0} (d\varphi) (x) \, ,}
$$
where $\NN_{1\to 0}$ is bounded from Radon measures to
functions of bounded variation. In other words,
$$
\widetilde \MM =\NN_{1\to 0} d \, .
$$
Similarly, we may decompose
$$
\widetilde {\MM^m} =\NN_{m,1\to 0} d \, ,
$$
and it is easy to see that
$$
\NN_{m,1\to 0} d \NN_{1\to 0}=\NN_{m+1,1\to 0} \, .
$$
(Just use that $d\sigma(x-\cdot)$ is the dirac at $x$.)
To finish,
$$
\widetilde {\MM^m} \widetilde \MM
=\NN_{m,1\to 0} d \NN_{1\to 0} d= \NN_{m+1,1\to 0}d= \widetilde {\MM^{m+1}}\, .
\hbox{\qed}
$$
\enddemo

\subhead 2.5 $\Det^\#(\Id-z\MM)$ and the spectrum of $\MM$
\endsubhead

\noindent
In this final section, we exploit the Milnor-Thurston identity to prove:

\proclaim{Theorem 6}  $\Det^\#(\Id-z\MM)$ is holomorphic in the
disc $|z| < 1/\widehat R$ and its zeroes in this disc are the inverses
of the eigenvalues of modulus larger than $\widehat R$ of $\MM$
acting on $\BB$. The order of the zero coincides with the algebraic
multiplicity of the eigenvalue.
\endproclaim

\demo{Proof of Theorem 6}
Combining Lemma 3 with the first equality in Theorem 4
$$\Det^\#(\Id -z \MM)=\Det_*(\Id +\widehat \DD(z))$$
and the fact that the spectral radius of $\widehat \MM$
on bounded functions is not larger than $\widehat R$, we get the holomorphy
claim ($\sigma(\cdot-y)$ is a bounded function for all $y$).

Using the second equality we see that if $z_0$ with $|z_0|< \widehat R$
is a zero of $\Det^\#(\Id -z \MM)=\Det_*(\Id +\DD(z))^{-1}$, then $1/z_0$
must be an eigenvalue of
$\MM$ acting on $\BB$ ($\sigma(\cdot-y)\in \BB$  for all $y$).

Let us now prove that if $\lambda_0=1/z_0$ is a simple eigenvalue,
then the order of the pole of $\Det_*(\Id +\DD(z_0))$ is at most one.
Writing $\PP:\BB \to \BB$ for the rank-one spectral projector
associated to $\MM$ and $\lambda_0$, we may decompose
$$
(\Id-z\MM)^{-1}=
{1\over 1-\lambda_0 z} \PP+  (\Id-z\MM)^{-1} (\Id-\PP) \, ,
$$
the second term being holomorphic in a neighbourhood of $z=z_0$.
We wish to use the above decomposition via multilinearity of the
determinants. For this it is useful
to use a Plemelj-Smithies formula
$$
\Det_*(\Id+\DD(z))=1+ \sum_{n=1}^\infty {1\over n!} \Phi_n(\DD(z)) \, ,
$$
where, writing $\KK^z(x,y)$ for the kernel of
$\DD(z)$, we set for $n\ge 1$
$$
\Phi_n(\DD(z)) =\int_{I^n} \det_{n\times n}
\biggl (\KK^z(x_i, x_j) \biggr ) \, d\mu(x_1) \cdots d\mu(x_n) \, .
$$
The Plemelj-Smithies formula claimed above can be obtained
by combining the consequence of Lemmas~5--6 and the Plemelj-Smithies
formula for
$\Det_2(\Id+\DD(z))$ (see Corollary of Theorem~2 in the Appendix).
Using our decomposition of the resolvent and the definition
of $\KK^z(x,y)$, we find
$$
\eqalign{
\KK^z(x,y)&=
\sum_\omega z {dg_\omega \over d\mu} (x) [(\Id-z\MM)^{-1} \sigma(\cdot-y)]
\psi_\omega(x)\cr
&= {z\over 1-\lambda_0 z} \alpha(x) \beta(y) + B^z(x,y) \, ,
}
$$
where $\alpha$ and $\beta$ are independent of $z$ and are
bounded on $I$, while $B^z(x,y)$ is bounded on
$I \times I$ and depends holomorphically
on $z$ in a neighbourhood of $z_0$.
We next develop each
$$
\det_{n\times n}
\biggl ({z\over 1-\lambda_0 z} \alpha(x_i) \beta(x_j)+\BB^z(x_i, x_j) \biggr )
$$
by multilinearity. The Hadamard inequality (used in the classical Fredholm
theory) gives
$$
\det_{n\times n}
\biggl (\BB^z(x_i, x_j) \biggr ) \le C_0^n n^{n/2} \, ,
$$
for some finite constant $C_0$, in a neighbourhood of $z_0$.
The other terms have one
or several columns of the form
${z\over 1-\lambda_0 z} \alpha(x_i) \beta(x_j)$. If there is a single
such column, the Hadamard inequality gives that in a neighbourhood
of $z_0$ the determinant
is at most
$$
{C_1 \over |1- \lambda_0 z|} C_0^n n^{n/2} \, ,
$$
where $C_1$ is another finite constant. If there are two or more such
columns, they are proportional so that the corresponding determinant vanishes.
Finally, we get by summing all terms and integrating over our finite
measure $\mu^n$:
$$
|\Phi_n(\DD(z))| \le { C_2^{n+1} \over |1-\lambda_0 z|} n^{n/2+1} \ .
$$
Putting this estimate back into the Plemelj-Smithies formula, we see that
the order of the pole at $z_0=1/\lambda_0$ is at most one, as claimed.
Note that the argument may be adapted if the algebraic multiplicity
is larger than one, showing that the order of the pole
is at most the algebraic multiplicity of the eigenvalue
(we shall not need this).

\medskip
To finish the proof, we show (using again the first equality
in Theorem~4) that if $\lambda_0$ with $|\lambda_0|> \widehat R$
is an eigenvalue of $\MM$ acting
on $\BB$, then $z_0=1/\lambda_0$ is a zero of
$\Det^\#(1-z\MM)=\Det_*(1+\widehat \DD(z))$ of order the algebraic multiplicity
of the eigenvalue. Since $\Det_*(1+\widehat \DD(z))=\Det_*(1+\widehat \DD(z)^*)$,
it is enough to show that $-1$ is an eigenvalue of $\widehat \DD(z_0)^*$
acting on $L^2(d\mu)$, with the correct multiplicity.

Let then $\varphi \in \BB$ be an eigenfunction for $\MM$ and the
eigenvalue $\lambda_0$. We can assume that $\varphi$ has only regular
discontinuities, and the eigenfunction equation implies that $\varphi$
is supported in $\cup_\omega I_\omega$. In particular,
$\varphi \in L^2(d\mu)$. We next show that $\widehat \DD(z_0)^*\varphi=-\varphi$.
First recall that
$$
\eqalign{
\widehat \DD(z_0)^*\varphi(y)&=
\int \sum_\omega z\epsilon_\omega {d(g_\omega \circ \psi_\omega^{-1})\over d\mu}(x)
[(\Id-z\widehat \MM)_x^{-1} \sigma(\cdot-y) ]\psi_\omega^{-1}(x)
\varphi(x) \, d\mu(x)\cr
&=
\int \sum_\omega z\epsilon_\omega
[(\Id-z\widehat \MM)_x^{-1} \sigma(\cdot-y) ]\psi_\omega^{-1}(x) \varphi(x)
d(g_\omega \circ \psi_\omega^{-1})(x)\, .\cr
}
$$
Using
$$
\varphi d(g_\omega \circ \psi_\omega^{-1})
=d(\varphi (g_\omega \circ \psi_\omega^{-1}))  -
(d\varphi) (g_\omega \circ \psi_\omega^{-1})\, ,
$$
and the fact that $\varphi$ only has regular discontinuities, we get
$$
\eqalign{
(1+\widehat \DD(z_0)^*)\varphi (y)&=-\int d\varphi(x) \biggl ( \sigma(x-y)\cr
&\qquad +\sum_\omega z_0 \epsilon_\omega d(g_\omega \circ \psi_\omega^{-1})(x)
[(\Id-z_0\widehat \MM)_x^{-1} \sigma(\cdot-y)]\psi_\omega^{-1}(x)\biggr ) \cr
&\qquad+\int\sum_\omega
d(z_0\epsilon_\omega \varphi (g_\omega \circ \psi_\omega^{-1})) (x)
[(\Id-z_0\widehat \MM)_x^{-1} \sigma(\cdot-y)]\psi_\omega^{-1}(x)\cr
&=\int -d\varphi(x)  (\Id + z_0 \widehat \MM_x (\Id-z_0 \widehat \MM)_x^{-1})
\sigma(x-y)\cr
&\quad +\int \sum_\omega d(z_0 (\varphi \circ \psi_\omega) g_\omega )(u)
(\Id-z_0 \widehat \MM)^{-1}_u
\sigma(u-y)\cr
&=-\int d(\varphi - z_0 \MM \varphi)(x) (\Id-z_0 \widehat \MM)^{-1}_x
\sigma(x-y)=0\, .
}
$$
If $\lambda_0$ is a simple eigenvalue then we are done.
(More generally, the above computation shows that the order of the zero
is at least the geometric multiplicity of the eigenvalue.
However, contrarily to the claims in the end of the proof
of Theorem 4.4.5 [Go], it is not clear how to relate
directly the order of the zero and the algebraic multiplicity.) Otherwise,
letting
$m_0$ be the algebraic multiplicity of $\lambda_0$, we  may find a small
perturbation $\MM_\delta$ of
$\MM$ (within the class of transfer operators associated to data $\psi_\omega$,
$g_\omega$) such that $\lambda_0$ is replaced by $m_0$
simple eigenvalues for $\MM_\delta$. (The details are left to the reader.
We use that a small perturbation in operator norm does not change the
spectrum away from a neighbourhood of $\lambda_0$ too much and cannot
increase the algebraic multiplicity of any perturbation
$\lambda_{0,\delta}$ of $\lambda_0$, and
also that the possibility of choosing independent supports for
the $\psi_{\omega_\delta,\delta}$
gives us enough ``degrees of freedom,'' to ensure that
$\Ker (\Id - \lambda \MM_\delta)^{m_0}$ is one-dimensional
for $\lambda$ in a neighbourhood of $\lambda_0$.)
Then the arguments already given show that $\Det^\#(1-z \MM_\delta)$
has $m_0$ simple zeroes in a neighbourhood of $z_0$, and that the
holomorphic  functions
$\Det^\#(1-z\MM_\delta)$ converge uniformly to $\Det^\#(1-z\MM)$
in compact sets.
\qed
\enddemo

\remark{Exercise 11}
In fact, any eigenfunction $\varphi$ in $\BB$ for $\MM$ and an eigenvalue
$\lambda$ of modulus larger than $\widehat R$ has a continuous representative.
(The proof of this fact is analogous to the proof of Lemma~4 in Chapter~1.
The starting point is to introduce $\tilde \varphi(x)=\varphi(x+)-\varphi(x-)$,
noting that $\MM \tilde \varphi=\lambda \varphi$, and writing
$(\Phi,\tilde \varphi)
=\sum_{x} \Phi(x)\tilde \varphi(x)$ for any bounded function $\Phi$.
Then, one can show that for each
bounded $\Phi$ we have $(\Phi, \tilde \varphi)
=\lambda^{-1}(\widehat \MM \Phi, \tilde \varphi)$ and end by iterating.)
\endremark

\remark{Exercise 12}
In Chapter 1 we used the decomposition
$$
d_{0\to1} \MM \SS_{1\to 0}= \widehat \MM_{1\to 1} +
\widehat \NN_{0\to 1} \SS_{1\to 0}
$$
to prove Theorem~1, with $\widehat \NN_{0\to 1}(\varphi)
=\sum_\omega dg_\omega \, \varphi \circ \psi_\omega$. We can also write
$$
\MM = \SS_{1\to 0} \widehat \MM_{1\to 1} d_{0\to1}+
\SS_{1\to 0}\widehat \NN_{0\to 1} \, ,
$$
with $\SS_{1\to 0}\widehat \NN_{0\to 1}$  compact on $\BB$.
Show that
$$
\SS_{1\to 0} \widehat \MM_{1\to 1} d_{0\to1}=\widehat \MM_{0\to 0}\, ,
$$
and (applying $d_{0\to 1}$ to both sides)  that
$$
\SS_{1\to 0}\widehat \NN_{0\to 1} =\SS_{0\to 0} \widehat \NN_{0\to 0}\, ,
$$
with $\widehat \NN_{0\to 0}\varphi =\sum_\omega \epsilon_\omega
{dg_\omega \circ \psi_\omega^{-1}\over d \mu}\varphi \circ \psi_\omega^{-1}$
and $\SS_{0\to 0}\varphi(x)=\int \sigma(x-y)\varphi(y)\, d\mu(y)$.
This implies that (in the notations of Chapter~2)
$$
(\Id -z \MM)^{-1}=
(\Id -z  \widehat \MM )^{-1}(\Id-z \SS\widehat \NN
(\Id -z  \widehat \MM )^{-1})^{-1}\ , .
$$
Finally, recall from the proof of Theorem~4 that
$$
\Det_* (\Id+\widehat \DD(z))=\Det_*
(\Id +z  \SS\widehat \NN  (\Id -z \widehat \MM )^{-1})\, .
$$
\endremark


\head 3. Kneading theory in higher dimensions
\endhead

In this last chapter, we shall discuss partial extensions of the
results of Chapter~ 2 to higher dimensions.
In fact we shall only present a higher-dimensional
version of the Milnor and Thurston identity
(Ph.D. of Baillif [Bai] based on an unpublished idea
of Kitaev), under a transversality assumption.

\subhead 3.1 Setting -- The higher-dimensional Milnor-Thurston formula
\endsubhead

In this chapter $\Omega$ is as before a
finite index-set, and $n \ge 2$ denotes the dimension,
i.e., we are going to work in a compact subset
$K$ of $\real^n$.
We also fix an integer order of differentiability
$r \ge 1$.
To each $\omega \in \Omega$ we associate an
nonempty open set $U_\omega \subset \real^n$ and
a (local) $C^r$ diffeomorphism
$$
\psi_\omega : U_\omega \to \psi_\omega(U_\omega) \, ,
$$
assuming that
$\bigcup_\omega U_\omega \cup \bigcup_\omega \psi_\omega(U_\omega)
\subset K$. We also consider a function $g_\omega :\real^n \to
\complex$ satisfying (these are not the weakest possible requirements)
\roster
\item
$g_\omega$ is supported in $U_\omega$,
\item
$g_\omega$ is $C^r$.
\endroster

Additionally, we make a {\it transversality} assumption
on the dynamics $\{ \psi_\omega\, , \mid \omega \in \Omega\}$:
For each $x \in \real^n$ such that there exist $m \ge 1$
and $\vec \omega=\omega_m \cdots \omega_1$ with
$\psi^m_{\vec \omega} (x)=
\psi_{\omega_m} \circ \cdots \circ \psi_{\omega_1}(x)=x$
(in particular, $x$ lies in the compact set $K$),
we have
$$
1 \notin \sp(D_x \psi_{\vec \omega}^m)\, .
$$
(In other words, $\psi^m_{\vec \omega}$ is a transversal
diffeomorphism.)

Note that in dimension one no such assumption was
present. We expect that this transversality requirement
(which is weaker than hyperbolicity) will be
eventually suppressed, but this will require some
additional work.

Here are two consequences of  transversality:
First, since the fixed points of a
transversal diffeomorphism are isolated,
and since $K$ is compact, there are only finitely
many fixed points of $\psi^m_{\vec \omega}$
for each fixed $m \ge 1$. Also, the Lefschetz number
$L(x, \psi_{\vec \omega}^m)$ does not vanish and can
be written
$$
L(x, \psi^m_{\vec \omega})= \sgn (\det (\Id - D_x \psi_{\vec \omega}^m))
\in \{ +1, -1\} \, .
$$
\smallskip

Before we introduce the kneading operators and state
the higher-dimensional version of the Milnor and Thurston
formula, we need to recall some notations and definitions.
Our $C^r$ assumption will allow us to use Lebesgue measure
$dx$ as a reference measure and $L^q$ will always denote
$L^q(dx)$.

We will be working not only with functions,
but more generally with $k$-forms.
For $k=0, \ldots, n$ and $1 < q < \infty$, we write
$A_k$ and $A_{k, K}$ for the vector spaces of $k$-forms
(with $C^\infty$ coefficients), respectively $k$-forms
supported in $K$. Also,
we will write $\AA_{k,L^q}$ and $\AA_{k, L^q(K)}$ for the
vector spaces of $k$-forms on $\real^n$ with $L^q$ coefficients,
respectively $L^q$ coefficients supported in $K$.
(We refer to [Sp] for the basic theory of differential
forms.)
This vector space inherits a Banach norm  from
the norms of the coefficient functions. Note however that the
corresponding Banach space is not convenient for
the spectral theory of the transfer operators $\MM_k$
to be introduced in a moment. It is however useful
for intermediate steps, in particular when considering
the kneading operators, also to be introduced below.

We shall denote the exterior derivative from
$A_{k, (K)}$ to $A_{k+1, (K)}$ by $d_k$ (or simply
$d$ when there is no ambiguity). Recall that if
$\phi = \sum_{\vec \jmath \in I(k)} \phi_{\vec \jmath}
\, dx_{\vec \jmath}$
(where $I(k)$ denotes the set of ordered
$k$-tuples in $\{1, \ldots, n\}$ and $dx_{\vec \jmath}=
dx_{j_1} \cdots dx_{j_k}$),
then $d\phi = \sum_{i=1}^n {\partial  \over \partial x_i}
\phi_{\vec \jmath} \, dx_{\vec \jmath} \wedge dx_i$ and that
$d_{k+1} d_k=0$. Sometimes $d_k$
will be considered on forms whose coefficients
are not $C^\infty$.  We shall work with the pull-back
$\psi_\omega^* $ on $\AA_k$ (or $\AA_{k,K}$,
or $\AA_{k,L^q}$, $\AA_{k,L^q(K)}$) of $\psi_\omega$.

We may now associate an $(n+1)$-tuple of transfer operators
$\MM=(\MM_k\, , k=0, \ldots n)$ to the data $g_\omega$, $\psi_\omega$,
where $\MM_k$ acts on $\AA_{k,L^q(K)}$  (for example)
by setting
$$
\MM_k \phi =
\sum_{\omega} g_\omega  (\psi_\omega^* \phi) \, .
$$
For $k=0$, we recover the previous definition:
$$
\MM_0 \phi(x) = \sum_\omega g_\omega(x) \, (\phi \circ \psi_\omega)(x) \, .
$$
For $m \ge 1$ we write $\MM^m$ for the $(n+1)$-tuple
$(\MM^m_k\, , k=0, \ldots n)$. Clearly $\MM^m$ is associated
to the data $\Omega^m$, $\psi_{\vec \omega}^m$
(when the domain of definition of the composition
is not empty) and
$$
g_{\omega_m}(\psi^{m-1}_{\vec \omega})
\cdots g_{\omega_2} (\psi_{\omega_1}) g_{\omega_1} \, .
$$
(Note that we do not claim, nor shall we need, that
the  $g_\omega$, $\psi_\omega$ are unambiguously
determined by the operators $\MM_k$.)

Let us now define the sharp trace and the sharp determinant:
\definition{Definition 1}
Let $\MM$ be associated to data $\{\psi_\omega\, ,
g_\omega, \omega \in \Omega \} $ as above,
then the sharp determinant of $\MM$ is defined by
$$
\Det^\# (1-z \MM) = \exp - \sum_{m=1}^\infty
{z^m \over m} \Tr\# \MM^m  \, ,
$$
where
$$
\Tr^\# \MM = \sum_{\omega\in \Omega} \sum_{x \in \Fix \psi_\omega}
g_\omega(x) L(x, \psi_\omega) \, .
$$
\enddefinition

The formula in the following exercise will  play a
part in the proof of the main theorem of the present section:

\remark{Exercise 1}
Define (here this is only a notation)
for each $k=0, \ldots, n$
$$
\Tr^\flat \MM_k= \sum_\omega \sum_{x \in \Fix \psi_\omega}
g_\omega(x) {\Tr \Lambda^k (D_x \psi_\omega)
\over |\det (\Id -D_x \psi_\omega)|} \, ,
$$
and
$$
\Det^\flat(\Id-z\MM_k)
=\exp - \sum_{m=0}^\infty{z^m\over m}
\Tr^\flat \MM^m_k \, .
$$
Show that
$$
\Det^\# (\Id- z\MM) = \prod_{k=0}^n (-1)^k
\Det^\flat(\Id-z\MM_k) \, .
$$
(Hint: use that for a finite matrix $A$
we have $\Det (\Id-A) = \sum_{k=0}^n (-1)^k \Tr\Lambda^k A$.)
\endremark

\smallskip

{\bf The Milnor-Thurston formula}
\smallskip

We are going to define homotopy operators
$$
\SS_k : \AA_{k+1, C^{r(-1)}(K)} \to \AA_{k, C^{r(-1)}} \, ,
k=-1, \ldots , n-1 \, ,
$$
(it is in fact possible to see that
$\SS_k(\AA_{k+1, L^q})\subset \AA_{k, L^q}$
for each $1 <q < \infty$) with
the property that, on compactly supported $k$-forms,
$$
d_{k-1}\SS_{k-1}+\SS_k d_k=\Id \, .
$$
The above homotopy equation can be solved because we are considering
forms in $\real^n$. In order to apply the techniques presented in this
chapter to dynamical systems on compact manifolds, one should
first embed the manifold in $\real^n$ for suitable $n$ and then
extend the dynamics in a tubular neighbourhood of the manifold.
(See [Bai].)
It will be clear from the construction below that the kernel $\sigma_k(x,y)$
of $\SS_k$ is smooth except on the diagonal $x=y$
in $\real^n$ where its singularities are of the type
$x_j/\|x\|^{n}$.

We shall also introduce auxiliary transfer-type operators
$$
\NN_k : \AA_{k, L^q}\to \AA_{k+1, L^q}\, , k=0, \ldots, n-1\, ,
1 < q < \infty \, ,
$$
defined by
$$
\NN_k \phi = (d_k \MM_k -\MM_{k+1} d_k) \phi
= \sum_\omega dg_\omega \wedge (\psi_\omega^* \phi)
$$
(we used the Leibniz formula). The operators $\NN_k$ also
map $\AA_{k, C^{r-1}}\to \AA_{k+1, C^{r-1}(K)}$.

Finally, the kneading operators are defined, for the moment
as formal power series with coefficients bounded operators
from $\AA_{k+1, C^{r-1}(K)}$ to $\AA_{k+1, C^{r-1}(K)}$ (with
$k=0, \ldots, n-1$) by
$$
\DD_k(z)=z \NN_k (\Id-z \MM_k)^{-1} \SS_k \, .
$$
Writing $\DD_k(z)$ as a kernel operator
with kernel $\KK^z_k(x,y)=\sum_{j=0}^\infty \kappa_{k,j}(x,y) z^j$
we shall prove (using the transversality assumption, see Lemma~4)
that $\kappa_{k,j}(x,x)\in L^1(\real^n)$ and define the formal trace
$\Tr_*(\DD_k(z))$ to be the  power series
$$
\Tr_*(\DD_k(z))=\sum_{j=0}^\infty z^j
(-1)^{(n+1)k} \int_{\real^n} \kappa_{k,j}(x,x) \, dx \, .
$$
(Note that the sign factor is
not present e.g. in odd dimensions.)
Proceeeding similarly for iterates of $\DD_k(z)$ we can
define a formal determinant from the formal trace as usual:
$$
\Det_*(\Id+\DD_k(z))=
\exp -\sum_{\ell=1}^\infty {z^\ell\over \ell}
\Tr_* (\DD_k(z)^\ell) \, .
$$
We shall then
prove the following theorem:

\proclaim{Theorem 1 (Milnor-Thurston-Kitaev-Baillif formula [Bai])}
In the sense of formal power series:
$$
\Det^\#(\Id-z\MM) =
\prod_{k=0}^{n-1} \Det_*(\Id+\DD_k(z))^{(-1)^{k+1}} \, .
$$
\endproclaim

\remark{Remarks on Theorem 1}
\roster
\item
We shall see that in fact for small enough $|z|$
the operator $\DD_k(z)$ is bounded on $\AA_{k+1, L^q}$
for $1 < q < \infty$ and that $\DD_k(z)^{[n/2]+1}$ is
Hilbert-Schmidt on $L^2(K)$. This additional information
allows us to express $\det_*(\Id+\DD_k(z))$ as the product
of a regularised determinant of order $2[n/2]+2$ and a holomorphic
non-vanishing function. This is useful to show that
$\Det_*(\Id+\DD_k(z))$ has a positive radius of convergence
(see [Bai])
and to extract spectral information from its zeroes [BB].
\item
By using the transversality assumption
it can be shown [Bai] that if the data $\psi_\omega$ and $g_\omega$
are $C^\infty$ then $\Det_*(\Id+\DD_k(z))$ is in fact the
flat-determinant [AB1, AB2] of $\DD_k(z)$. Properties of
the flat determinants can be used to give a short proof  of Lemma~6 below.
One has to use an approximation argument in case the
original data is just $C^r$ for finite $r$ (see [Bai]). 
We shall thus use from now on the notation (which is slightly
abusive if $r \ne \infty$):
$$
\Det^\flat(\Id+\DD_k(z)) = \Det_*(\Id+\DD_k(z)) \, , 
\quad \Tr^\flat \DD_k(z)=\Tr_* \DD_k(z) \, .
$$
\endroster
\endremark

\smallskip
\remark{Exercise 2} For $n \ge 2$, show that the functions
$x_j/\|x\|^n$ are in $L^q(K')$ for any compact subset $K'$
of $\real^n$ and all $1 \le q < n/(n-1)$.
\endremark

\subhead 3.2 Definition of the homotopy operators $\SS_k$
\endsubhead

Let us now proceed with the definition of the homotopy operators
$\SS_k$. The starting point is to find an inverse to the
Laplacian acting on (compactly supported) $k$-forms i.e. a  solution $G_k$
(for $k=0,\ldots, n$) to
$$
\Delta G_k = \Id \, , \quad G_k \Delta = \Id
$$
on $\AA_{k,C^\infty(K)}$, where $\Delta=\Delta_k$ is the
Laplacian operator acting on $k$-forms:
$$
\Delta_k \biggl ( \sum_{\vec \jmath \in I(k)}
 \phi_{\vec \jmath} dx_{\vec \jmath}\biggr )
= - \sum_{\vec \jmath \in I(k)} \biggl ( \sum_{i=1}^n {\partial^2 \over \partial x_i^2}
\phi_{\vec \jmath}\biggr ) \, dx_{\vec \jmath}\, .
$$

\proclaim{Lemma 2} Let $E \in \AA_{n, L^1(K)}$ be the ``Green kernel''
$$
E(x)=e(x) dx_1 \wedge \cdots dx_n
= \cases
{\Gamma (n/2) \over (n-2) 2 \pi^{n/2}}
{1\over \|x\|^{n-2}} \, dx_1 \wedge \cdots dx_n& n \ge 3\cr
{1 \over 2\pi} \log (\|x\|) \, dx_1 \wedge dx_2 &n=2 \, ,
\endcases
$$
where $\Gamma$ is Euler's gamma-function. For $k=0, \, \ldots, \, n$
define a $k$-form in $x$ and an $n-k$-form in $y$
(with coefficients in $L^1(K)$)  $E_k(x,y)$ by
$$
E(x-y) = \sum_{k=0}^n (-1)^{n(k+1)} E_k(x,y) \, .
$$
Then the operator on compactly
supported $k$-forms with $C^{t}$ coefficients
($t \ge 0$) defined by
$$
G_k \phi (x) = \int_{\real^n } E_k(x,y) \wedge \phi(y)
$$
is an inverse for the Laplacian $\Delta_k$.
\endproclaim

\demo{Proof of Lemma 2}
The function $e(x)$ in the Green's kernel has the property that
(as a distribution on compactly supported $C^\infty$ functions)
$$
\Delta e(x) = \delta_0 \, ,
$$
where the right-hand-side is the dirac mass at $0$.
This property can be proved by using  Green's formula --
see [Sch, p. 46] for details of this classical and elementary computation.
From this, it is not difficult to deduce that $\Delta G_k =\Id$
by noting first that
$$
E_k(x,y) =\sum_{\vec \jmath \in I(k)}
s(\vec \jmath ') e(x-y) \, dx_{\vec \jmath} \wedge dy_{\vec \jmath'}
$$
where $\vec \jmath'\subset I(k)$ is the ordered
complementary of $\vec \jmath$ in
$\{ 1, \ldots , n\}$ and $s(\vec \jmath')\in \{ -1, +1\}$
is the sign of the permutation reordering
$(\vec \jmath', \vec \jmath)$;
so that
$$
\Delta_x E_k(x,y)=
\delta_{0,x} (x-y)
\biggl ( \sum_{\vec \jmath } s(\vec \jmath')
dx_{\vec \jmath} \wedge dy_{\vec \jmath'} \biggr ) \, .
$$
Indeed, it follows
that for any $k$-form $\phi=\phi_{\vec \ell} dx_{\vec \ell}$
$$
\eqalign{
\Delta G_k \phi (x)&= \int_{\real^n}\Delta_x E_k(x,y) \wedge \phi(y)\cr
&=  \int_{\real^n}\Delta_y E_k(x,y) \wedge \phi_{\vec \ell}(y)
\wedge dy_{\vec \ell} \cr
&= {s(\vec \ell')}
\int_{\real^n} \delta_{0,y} (x-y)  \phi_{\vec \ell}(y) \,
dx_{\vec \ell}  \wedge dy_{\vec \ell'}  \wedge dy_{\vec \ell} \, .
}
$$
Integration by parts and one more use of
${\partial^2\over dy_i^2}E(x-y)={\partial^2\over dx_i^2}E(x-y)$
then implies $G_k \Delta =\Id$.
\qed \enddemo

\smallskip

Recall now the classical identity
$$
\Delta=\Delta_k = d^*_{k+1} d_k + d_{k-1} d^*_k \, ,
$$
where $d^*_{k+1}: \AA_{k+1} \to \AA_k$ may be defined by duality
$$
<d^*_{k+1} \phi, \psi > = <\phi, d_k\psi>\, ,
$$
where, for any two $\ell$-forms $\varphi_1, \varphi_2$ we set
$$
<\varphi_1 , \varphi_2>=
\cases
\int \varphi_{1, \vec \jmath} \varphi_{2, \vec \jmath}
dx_1 \wedge dx_n  &\text{if } \varphi_1, \, \varphi_2
\text{ have the same support } \vec \jmath \in I(\ell)\cr
0&\text{otherwise} \, .
\endcases
$$
Note that if $\phi$ is a $C^1$ function then
$$
d^*_k \phi(x) \, dx_1 \wedge \cdots \wedge dx_k
= \sum_{j=1}^k  (-1)^{j+1} {\partial \over \partial x_j}
\phi(x) \, dx_1 \wedge \cdots
\wedge \widehat{dx_j} \cdots \wedge dx_k \, ,
$$
where $\widehat{dx_j}$ means that the factor $\widehat{dx_j}$
has been suppressed.

Use of the following homotopy operators was first suggested by Kitaev,
the expression given in the definition below was remarked by Ruelle
but the operators are the same as those appearing in [Bai]:

\definition{Homotopy operators}
For $k=-1, \ldots , n$ (and $t\ge 0$) we set
$$
\SS_k = d^*_{k+1} G_{k+1} : \AA_{k+1, C^{t}(K)}
\to \AA_{k, C^t} \, .
$$
\enddefinition

\remark{Exercise 2}
Show (formally) that $\SS_{k-1} \SS_k\equiv 0$.
\endremark

\proclaim{Lemma 3}
\roster
\item
The homotopy operator $\SS_k$
on $\AA_{k, C^{r(-1)}(K)}$  admits an expression
in kernel form as
$$
\SS_k \phi(x)= \int \sigma_k(x,y) \wedge \phi(y)
$$
where $\sigma_k(x,y)$ is a $k$-form in $x$ and an
$n-k$ form in $y$ obtained from the $n-1$-form
$$
\sigma(x):=d^* E (x)=
\sum_{i=1}^n (-1)^{i+1} {x_i \over \|x\|^n }
dx_1 \wedge \cdots \wedge \widehat {dx_i }\wedge \cdots dx_n \, ,
$$
by using the decomposition
$$
\sigma(x-y) = \sum_{k=0} ^{n} (-1)^{nk} \sigma_k(x,y)
$$
\item
$d_{k-1} \SS_{k-1} + \SS_k d_k =\Id$ on
$\AA_{k, C^{r(-1)}(K)}$.
\endroster
\endproclaim

\demo{Proof of Lemma 3}
In view of the definitions,
the proof of \therosteritem{1} consists in checking that the signs
match, and this is left as an exercise to the reader. Let us prove
\therosteritem{2}, i.e. verify that
$
d d^* G_{k+1} + d^* G_{k+2}d =\Id
$.
But this is an easy consequence of the following identity
$$
d d^* G_{k+1} + d^* G_{k+2} d=\Delta G_{k+1}
- d^* (d G_{k+1} -G_{k+2} d) \, ,
$$
since $d G_{k+1} -G_{k+2} d=0$ (integrating by parts).
\qed
\enddemo

\subhead 3.3 Properties of the kneading operators and other
kernel operators
\endsubhead

In order to prove Theorem 1, we shall make use of the
transversality assumption to prove that the kneading
operators, and also some other related operators, are
such that either their kernel can be integrated along
the diagonal in the sense of an $L^1$ function (Lemma~4), or
(Lemma~5) that their generalised (Schwartz) kernel (which a priori
is only a current over $\real^{2n}$) can be restricted
to the ($n$-dimensional) diagonal where it gives rise to a distribution, which
can then be evaluated over the constant function $1$ (say).
In fact, it is convenient for part of the computations to
assume that $r=\infty$. If the original data only enjoys
finite smoothness, an approximation argument can be used
(thanks to transversality). We refer to [Bai] for this,
and will only present the proof of Theorem~1 in the case
$r=\infty$.

To proceed, we introduce two vector spaces of operators corresponding
to the two cases just discussed.
The definitions will ensure that $\DD_{k}(z)$ is a power series
with coefficients operators in the first space $\KK_{k+1}$. In the case $r=\infty$, all auxiliary
operators which will be introduced in the proof of Theorem~1
will be power series with coefficients operators in
the second space $\KK^d_k$. 

\definition{Definition (The  spaces $\KK_k$ and $\KK_k^d$)}
Let $\psi_\omega$, $g_\omega$ be as in Section~3.1
for some $r\ge 1$, let $\MM_j$, $\NN_\ell$ and $\SS_m$
be the operators defined above (acting on locally supported
forms with coefficients in $C^{r-1}(K)$). We say that
a (finite) composition of operators $\MM_j$, $\NN_\ell$ and $\SS_m$
is admissible if the degrees of the forms match.
Fix an integer $0 \le k \le n$.
\roster
\item
We write $\KK_k$ for the vector space of bounded operators
from $\AA_{k, C^{r-1}(K)}$ to $\AA_{k, C^{r-1}}$
generated by admissible compositions of $\MM_j$, $\NN_\ell$ and $\SS_m$,
with at least one $\SS$ factor, no two immediately successive $\SS$ factors, and  the
first or last factor of type  $\MM$ or $\NN$.
\item
If $r=r-1=\infty$, we write $\KK^d_k$ for the vector space of bounded operators
from $\AA_{k, C^{\infty}(K)}$ to $\AA_{k, C^{\infty}}$
generated by admissible compositions of $\MM_j$, $\NN_\ell$, $\SS_m$,
and $d_q$
with at least one $\SS$ factor, no two immediately successive $\SS$ factors, at least one $\MM_j$
or $\NN_\ell$ between two $d$s,
and the first or last factor of type $\MM$ or $\NN$.
\endroster
\enddefinition

\remark{Exercise 3 ($\KK_k$ and $\KK^d_k$)}
Check that $\DD_k(z)\in \KK_{k+1}[[z]]$ and that
$\MM_k\in \KK_k^d$. (Hint: use $d\SS+\SS d=\Id$.)
\endremark

\proclaim{Lemma 4}
If $\QQ \in \KK_k$ then $\QQ$ is a linear combination of kernel operators
$\QQ_i$ with
$$
\eqalign{
\QQ_i: \AA_{k, C^{r-1}(K)} &\to \AA_{k, C^{r-1}(K)}\cr
\QQ_i \varphi(x) &= \int h(x) K(x,y) \tilde h(y) \wedge \phi(y)\cr
}
$$
where
\roster
\item
there are $\hat s(i)\ge 1$ and $\Psi=\psi^{\hat s}_{\vec \omega}$ so that
$K(x,y)$ is a $k$-form in $x$ and an $n-k$-form in $y$,
which is $C^r$ except on $\Psi(x)=y$;
\item
$h$ and $\tilde h$ are $C^{r-1}$ functions on $\real^n$,
$h$ is supported in $K$ if the leftwards factor in $\QQ_i$
is not $\SS_k$  while $\tilde h$ is supported
in $K$ if the rightwards factor is not $\SS_{k-1}$;
\item
$\chi_K K(x,x) \in L^p(\real^n)$  for all $1\le p \le n/(n-1)$.
\endroster
\endproclaim

Lemmma~4 allows us to define the flat-trace of an element of $\KK_k$
or of $\KK_k[[z]]$ (in particular, $\Tr^\flat(\DD(z))$ is now well-defined):

\definition{Flat trace of kernel operators}
If $\QQ \in \KK_k$ then, using the notation from Lemma~4, we define $\Tr^\flat \QQ\in \complex$
by 
$$
\Tr^\flat \QQ=(-1)^{(n+1)k}\int_{\real^n} h(x) K(x,x) \tilde h(x)\,  dx\, .
$$
If $\QQ(z)=\sum_{j=0}^\infty z^j \QQ_j \in \KK_k[[z]]$ (convergent or not), then 
we set $\Tr^\flat \QQ(z)\in \complex[[z]]$
$$
\Tr^\flat \QQ(z)=
\sum_{j=0}^\infty z^j (-1)^{(n+1)k}\int_{\real^n} h_j(x) K_j(x,x) \tilde h_j(x)\,  dx\, .
$$
\enddefinition

\remark{Remark on the flat trace on $K_d$}
Let $\QQ \in \KK_k$. We shall not need the following facts:
\roster
\item
If $r=\infty$, it is possible to show that $\Tr^\flat\QQ$ coincides
with the classical Atiyah-Bott flat trace [AB1,AB2]. See [Bai].
\item
One can prove that $\Tr^\flat\QQ$ only depends on the values
of $\QQ$ on $C^\infty$ locally supported  $k$-forms.
\endroster
\endremark

\demo{Sketch of proof of Lemma 4}
The first two claims can be proved by induction on the number of
factors, with $s(i)$ being the number of non-$\SS$ factors.
(It is convenient in the proof to introduce a unified notation for
the operators $\MM$ and $\NN$ by writing
$$T\varphi(x)=
\sum_\omega \eta_\omega (x) \wedge \psi_\omega^* \varphi(x)
$$
with $\eta_\omega$ in $\AA_{\ell-k, C^{r-1}(K)}$ for $\ell=k$
or $k+1$.)

We concentrate on the proof of \therosteritem{3}. Our starting
point is the following easily proved  expression for $K(x,y)$ (use Lemma~3):
$$
K(x,y)=\int_{(\real^n)^s}
G(x, x^{(1)}, \ldots, x^{(s)}, y) \,
H(x, x^{(1)}, \ldots, x^{(s)}, y) dx^{(1)} \wedge \cdots \wedge dx^{(s)}\, ,
$$
where $s\le \hat s(i)$, $G(x, x^{(1)}, \ldots, x^{(s)}, y)$ is
$C^{r-1}$ and compactly supported on $\{ x\} \times (\real^n)^s
\times \{y\}$, while, setting $x^{(0)}=x$, $x^{(s+1)}=y$,
$H(x^{(0)}, x^{(1)}, \ldots, x^{(s)}, x^{(s+1)})$
can be written as a linear combination of expressions
$$
H_{\vec \jmath} (x^{(0)}, x^{(1)}, \ldots, x^{(s)}, x^{(s+1)})=
\prod_{t=0}^s
{\psi_t (x^{(t)})_{j_t} - x^{(t+1)}_{j_t} \over
\|\psi_t (x^{(t)}) - x^{(t+1)} \|^n}
$$
for suitable $1 \le j_1, \ldots, j_s \le n$. (Each $\psi_t$  is
a composition of finitely many $\psi_\omega$s.) 
It thus suffices to show that each
$$
|\prod_{t=0}^{s} \chi_J(x^{(t)})
H_{\vec \jmath} (x^{(0)}, x^{(1)}, \ldots, x^{(s)}, x^{(0)})|^p
$$
belongs to $L^1(\real^n)^{(s+1)}$.
The singularities of $H_{\vec \jmath}$ are isolated (by transversality),
there are thus finitely many of them in a compact set.
We shall content ourselves with proving local integrability
in the neighbourhood of the ``worse'' possible singularities
$\hat x$, i.e.
$$
\cases
\psi_t(\hat x^{(t)}) = \hat x^{(t+1)}\, , \forall t=0, \ldots, s-1\, ,
\cr
\psi_s(\hat x^{(s)})= \hat x^{(0)}\, .
\endcases
$$
(The task of checking that the singularities corresponding to the vanishing 
of some,
but not all,
of the $s+1$ factors in the denominator of $H_{\vec \jmath}$ 
are also locally integrable is left to the reader.)
Let us perform the change of variables
$$
\cases
w^{(t)}=
\psi_t( x^{(t)}) -  x^{(t+1)}\, , \forall t=0, \ldots, s-1\, ,
\cr
w^{(s)}= \psi_s( x^{(s)})-  x^{(0)}\, .
\endcases
$$
We shall check later (using transversality) that the Jacobian  $J(w)=
|\det ({d\over dx} w(x))|$ of the
above change of variables (which is obviously $C^{r-1}$) does
not vanish at $\hat w=w(\hat x)=0$.
In the new coordinates, we have (with $\delta(\epsilon)\to 0$
as $\epsilon \to 0$)
$$
\eqalign
{
\int_{\| x-\hat x\| \le \epsilon} &| H(x^{(0)}, \ldots, x^{(s)},x^{(0)})|^p \,
dx^{(0)}  \cdots dx^{(s)} \cr
&\le \int_{w \in \real^{n(s+1)}\, , \| \hat w\| \le \delta} {1\over J(w)}
\biggl |  \prod_{t=0}^s
{w^{(t)}_{j_t}\over \|w^{(t)}\|^n } \biggr |^p dw^{(0)} \cdots dw^{(s)}\cr
&\le C \prod_{t=0}^s \int_{y\in \real^n \, , \|y\|\le \delta}
\left ({|y_{j_t}| \over \|y\|^n}\right )^p \, dy <\infty\, ,
}
$$
since $p<n/(n-1)$. 

It remains to check that $J(0)\ne 0$. For this, we observe that $J=\det D$ with
$D(w)$ the $n\times n$ matrix with entries
$$
D(w)_{tu}=
\cases
D(\psi_t)_{x(w)} & t=u\, , \cr
-1 &u=t+1 \le s \hbox{ or } t=s \, , u=0 \, , \cr
0& \hbox{ otherwise} \, .
\endcases
$$
If $\det D(0)=0$, then there would exist a nonzero vector $v_t$, $t=0,\ldots, s$,
with $Dv=0$, i.e., $v_{t+1}=D(\psi_t)_{\hat x_t} v_t$ and $v_0=D(\psi_s)_{\hat x_s} v_s$.  But then,
$v_0=D(\psi^{s+1}_{\vec t})_{\hat x} v_0$, contradicting transversality at the periodic
point $\hat x$.
\qed
\enddemo

We shall not give the proof of the following claim, referring instead to
[Bai]. It relies on transversality.
([Bai]  uses of results of Guillemin and Sternberg
and the wave-front-set, and he notes that
the flat-trace in $\KK^d_k$ coincides with that of Atiyah-Bott [AB1, AB2].)

\proclaim{Lemma-Definition 5}
Let $\QQ\in \KK_k^d$ and let $K_\QQ(x,y)$ be its Schwartz kernel [Sch], which
is a $k$-current in $x$, and an $n-k$-current in $y$, with coefficients distributions
of finite order. Then $\delta(x-y) K_\QQ(x,y)$ is a compactly supported distribution
on $\real^{2n}$. It can thus be evaluated on the constant function $1$,
giving a meaning to the following definition:
$$
\Tr^\flat \QQ = (-1)^{(n+1)k} \int_{\real^n} K_\QQ(x,x) \, dx\, .
$$
\endproclaim

\remark{Exercise 4}
Give an expression for the Schwartz kernel of $\MM_k$.
Check that the flat trace of $\MM_k$ as defined in
Lemma-Definition 5 coincides with the formal definition of Exercise~1.
\endremark

The proof of Lemma~5 in [Bai] shows that it is legitimate to invoke
the Fubini theorem when manipulating the Schwartz kernels of elements
of $\KK^d_k$. (This is not obvious since these are not functional
kernels.) As a consequence, he  proves:

\proclaim{Lemma 6}
If $\QQ_1$ and $\QQ_2$
are finite compositions of $\SS_m$,
$\NN_p$ and $\MM_q$ so that  $\QQ_1 \QQ_2 \in \KK_k^d$ and $\QQ_2\QQ_1 \in \KK_\ell^d$, then
$$
\Tr^\flat \QQ_1 \QQ_2 =\Tr^\flat \QQ_2 \QQ_1 \, .
$$
\endproclaim

As an immediate consequence, we get

\proclaim{Corollary of Lemma 6} Under the assumptions of
Lemma~6:
$$
\Det^\flat(\Id-z\QQ_1 \QQ_2)=\Det^\flat(\Id-z\QQ_2 \QQ_1)\, .
$$
If, additionally, $k=\ell$ and $\QQ_1$, $\QQ_2\in \KK_k^d$ then:
$$
\Det^\flat(\Id-z\QQ_1-z\QQ_2+z^2 \QQ_1\QQ_2)=\Det^\flat(\Id-z\QQ_1) \, \Det^\flat(\Id-z \QQ_2)\, .
$$
\endproclaim

\remark{Exercise 5}
Formulate Lemma 6 and its Corollary for elements of
$\KK_k^d[[z]]$, such as $(1-z\MM_k)^{-1}$.
\endremark

\subhead 3.4 Proof of the Milnor--Thurston formula in the $C^\infty$
case
\endsubhead

Let us exploit Lemmas~4--6 from \S~ 3.3 to sketch a proof of Theorem~1 under the
additional assumption that $r=\infty$. (We refer to [Bai] for the general case
which uses an approximation argument due to Kaloshin.)
We start by rewriting $\Det^\flat(\Id+\DD_k(z))$:
$$
\eqalign{
\Det^\flat(\Id+\DD_k(z))&=\Det^\flat(\Id+z\NN_k(\Id-\MM_k)^{-1} \SS_k)\cr
&=\Det^\flat(\Id+z\SS_k\NN_k(\Id-\MM_k)^{-1} )\cr
&=\Det^\flat(\Id-z(\MM_k-\SS_k\NN_k ) )\, \Det^\flat((\Id-\MM_k)^{-1} )\cr
&= \Det^\flat(\Id-z(\MM_k-\SS_k\NN_k )) \, (\Det^\flat(\Id-\MM_k))^{-1} \, .\cr
}
$$
By Exercise~1, it thus suffices to check that
$$
\prod_{k=0}^{n-1}\Det^\flat(\Id-z(\MM_k-\SS_k\NN_k ))^{(-1)^k}
=\Det^\flat(\Id-z\MM_n)^{(-1)^{n-1}}\, .
$$
But this follows from
$$
\eqalign
{
\Det^\flat(\Id-z(\MM_k-\SS_k\NN_k ))&=\Det^\flat(\Id-z(\MM_k-\SS_kd \MM_k +\SS_k\MM_{k+1}d))\cr
&=\Det^\flat(\Id-z(d \SS_{k-1}\MM_k +\SS_k\MM_{k+1}d))\cr
&=\Det^\flat(\Id-z d \SS_{k-1}\MM_k)\, \Det^\flat(\Id -z\SS_k\MM_{k+1}d)\cr
&=\Det^\flat(\Id-z d \SS_{k-1}\MM_k)\, \Det^\flat(\Id -zd \SS_k\MM_{k+1})
}
$$
(in the third line we used $d^2=0$).  Indeed, it is clear that the factors in the
alternated product cancel, except for
$$
\Det^\flat(\Id-z d \SS_{-1}\MM_0)=1 \hbox{ and } \Det^\flat(\Id-z d \SS_{n-1}\MM_n)=
\Det^\flat(\Id-z \MM_n)\, .\hbox{\qed}
$$

\medskip

\heading Appendix
\endheading

\subhead A. Hilbert-Schmidt operators and their regularised determinants
\endsubhead

Let $H$ be a separable Hilbert space. We recall here the results that we need,
referring to [GGK] for proofs and for more statements.

\definition{Definition (Hilbert-Schmidt operator)}
A compact linear operator $A :H \to H$ is called Hilbert-Schmidt,
noted $A\in \SS_2$ if $B=A^* A$ is a trace-class operator, noted
$B \in \SS_1$. A compact linear operator $B$ on $H$ is called trace-class
if $\sum_{j=1}^\infty s_j(B)<\infty$, where the $j$th
singular number of $B$ is defined by $s_j(B):=\sqrt{\lambda_j(B^*B)}$,
with
$$
\lambda_1(B^*B)\ge \lambda_2(B^*B) \ge \cdots \ge \lambda_j(B^*B)\ge \cdots
> 0 \, ,
$$
the sequence of nonzero eigenvalues of $B^*B$, repeated according
to multiplicity.
\enddefinition

\remark{Equivalent definition}
A compact linear operator $A$ on $H$ is Hilbert-Schmidt if and only if
there is an orthonormal basis $\{ \varphi_j\}$ of $H$ for which
$\sum_j \| A \varphi_j \|^2 < \infty$. (The sum then converges for every
orthonormal basis of $H$.)
\endremark

We refer e.g. to [DS1] for a proof of the very classical result:
\proclaim{Proposition (Hilbert-Schmidt operators on $L^2(d\mu)$)}
Let $\mu$ be a nonnegative measure on a $\sigma$-algebra of a set $I$.
Let $K(x,y)$ be a measurable function on $I\times I$. Then
the (kernel) operator
on the Hilbert space $H=L^2(d\mu)$ associated to $K$ by
$$
A\varphi(x)=\int_I K(x,y)\, \varphi(y) \, d\mu(y)\,
$$
is Hilbert-Schmidt if and only if $K\in L^2(d\mu\times d\mu)$, i.e.
$$
\int_{I\times I} |K(x,y)|^2 \, d\mu(x)\, d\mu(y) < \infty \, .
$$
\endproclaim

We now return to our abstract  separable Hilbert space and discuss
traces and determinants. Norms on the so-called Schatten classes $\SS_1$
and $\SS_2$ are introduced in the following exercise:

\remark{Exercise 0} Let $H$ be a separable Hilbert space and write
$L(H)$ for the algebra of bounded linear operators on $H$.
Show that the expressions
$$
\|A \|_2:= \sqrt{\sum_j s_j^2(A)  } \, , \quad \|B \|_1:=\sum_j s_j(B)\,
$$
define norms on $\SS_2$, respectively
$\SS_1$, that $\SS_1$ is a complete subalgebra of
$L(H)$ for this norm:
$$
\| B B'\| _1 \le \|B\|_1 \|B'\|_1 \, ,
$$
and that $\SS_1 \subset \SS_2$ continuously.
\endremark

\proclaim{Lemma 0 ($\SS_2 \SS_2 \subset \SS_1$)}
If $A, A'$ belong to $\SS_2$ then $AA' \in \SS_1$ and
$$
\| AA'\|_1 \le \|A \|_2 \|A'\|_2 \, .
$$
\endproclaim

\demo{Sketch of proof of the lemma}
For each $k\ge 1$ one can easily show that
$$
\sum_{j=1}^k s_j(A A') \le \sum_{j=1}^k s_j(A)s_j(A') \, ,
$$
to finish, one applies the Cauchy-Schwarz inequality.
\enddemo

\smallskip

{\bf The algebra of trace-class operators $\SS_1$}

We already noted in Exercise 0 that $\SS_1$ is a subalgebra of
$L(H)$. This algebra is in fact continuously embedded in $L(H)$,
i.e. for each $B \in \SS_1$ the operator norm is bounded by the
norm in $\SS_1$:
$$\|B\|_{L(H)} \le \| B\|_1\, . $$
This embedded subalgebra has the approximation property that
the space of finite rank operators $\FF$ on $H$ is dense in $\SS_1$
(for the $\SS_1$ norm $\| \cdot \|_1$). We are thus in a position
to apply the following extension theorem (see e.g.  [GGK, Chapter II.2]
for a proof) to $\widetilde \SS=\SS_1$:

\proclaim{Theorem 1 (Extending the trace and determinant)}
Let $\widetilde \SS$ be a continuously embedded subalgebra
of $L(H)$ with the approximation property. The following properties
are equivalent:
\roster
\item
The trace $F \mapsto \Tr F$ on $\FF\cap \widetilde \SS$ is a bounded functional for
the norm $\|\cdot \|_{\widetilde \SS}$.
\item
The trace $\Tr F$ and the determinant $F \mapsto \Det (\Id+F)$ admit continuous
extensions from $\FF \cap \widetilde \SS$ to $\widetilde \SS$. The continuity is
in the sense of the $\widetilde \SS$ norm and we have for any sequence
$F_n \in \FF \cap \widetilde \SS$ converging to
$A \in \widetilde \SS$:
$$
\Tr(A)=\lim_{n \to \infty} \Tr F_n \, , \quad
\Det(\Id + A)=\lim_{n \to \infty} \Det(\Id + F_n)\, .
$$
\endroster
\endproclaim

We shall make use of the following properties of the extended
determinants (the proofs are to be found in [GGK, II.3 and II.6]):
\proclaim{Theorem 2 (Properties of the extended determinant)}
Assume that we are in the equivalent conditions of the previous theorem. Then
for each compact  $A \in \widetilde \SS$:
\roster
\item
The function $\lambda \mapsto \Det(\Id-\lambda A)$ is an entire
function. (There is a formula for the coefficients of the Taylor series
at zero, called the Plemelj-Smithies formula.)
\item
$\Det(\Id-\lambda A) = \exp - \sum_{n=1}^\infty {\lambda^n \over n}
\Tr A^n \, .$
\item
$\Det(\Id -\lambda_0 A)=0$ with order $m_0 \ge 1$ if and only if
$1/\lambda_0$ is an eigenvalue of $A$ of algebraic multiplicity $m_0$.
\endroster
\endproclaim

The following property is essential to our application:
\remark{Exercise 1 (Analyticity of the extended trace and determinant)}
Let $z \mapsto A(z)$ be an analytic map at $z_0 \in \complex$
with each  $A(z)$ in a Banach algebra $\widetilde \SS$ satisfying the
equivalent  conditions of the
extension Theorem~1. Then both maps $z \mapsto \Tr A(z)$ and
$z \mapsto \Det(\Id+A(z))$  are analytic at $z_0$.
\endremark

\smallskip

Finally, we have:

\proclaim{Lidskii Trace Theorem ([GGK, IV.6])}
For $A\in \SS_1$, writing $\lambda_j(A)$ for the eigenvalues
of $A$  repeated with multiplicity, we have
$$
\Tr A= \sum_j \lambda_j(A)\, , \quad
\Det (\Id -A) = \prod_j (1-\lambda_j(A)) \, .
$$
\endproclaim

\smallskip

{\bf Hilbert-Schmidt operators and their regularised determinants}

If $A$ is a Hilbert-Schmidt operator on a separable Hilbert space $H$,
then the following operator is trace-class:
$$
R_A:= \Id - (\Id-A) \exp(A)\, .
$$
Indeed
$$
\eqalign{
R_A &= \Id -  \sum_{j=0}^\infty{ A^j\over j!}+
\sum_{j=0}^\infty{ A^{j+1}\over j!}
\cr
&=\sum_{j=2}^\infty{ A^j (1-j) \over j!} \, ,}
$$
so that $R_A$ is an absolutely convergent sum of operators
in $\SS_1$ (use Lemma~1).

\definition{Definition (regularised determinant)}
To $A\in \SS_2$ we associate a regularised determinant (of order two) by
setting:
$$
\Det_2(\Id-A) = \Det (\Id-R_A)=\Det((\Id-A)\exp(A))\, .
$$
\enddefinition

Note that there exists a theory of regularised determinants
of order $p\ge 2$ for the Schatten classes $\SS_p$ which have the
property that $A\in \SS_p$ implies $\AA^p \in \SS_1$.
(We refer to [GGK].)

The regularised determinant immediately inherits several properties from
the determinant in $\SS_1$:

\proclaim{Corollary of Theorem 2 (Properties of the regularised
determinant)}
For each $A \in  \SS_2$:
\roster
\item
If $A\in \SS_1$ (for example $A \in \FF$) then
$\Det_2(\Id-A)=\Det(\Id-A) \cdot \exp(\Tr A)$.
\item
The function $\lambda \mapsto \Det_2(\Id-\lambda A)$ is an entire
function. (There is a formula for the coefficients of the Taylor series
at zero, also called the Plemelj-Smithies formula [GGK,IX.3].)
\item
$\Det_2(\Id-\lambda A) = \exp - \sum_{n=2}^\infty {\lambda^n \over n}
\Tr A^n \, .$ [GGK, IX.3]
\item
$\Det_2(\Id -\lambda_0 A)=0$ with order $m_0 \ge 1$ if and only if
$1/\lambda_0$ is an eigenvalue of $A$ of algebraic multiplicity $m_0$.
\endroster
\endproclaim

As a consequence e.g. of \therosteritem{1} we see that the regularised
determinant is {\it not multiplicative.} In applications, it is often
necessary to complete it by a factor ``replacing'' the missing $\exp-\Tr$.

\demo{Proof of the Corollary}
We only prove \therosteritem{4}, leaving the other claims as exercises.
(In particular, \therosteritem{3} follows from \therosteritem{1}
and the fact that $\Det_2(\Id-A)=\lim \Det_2(\Id-F_n)$
where $F_n$ is a sequence of finite-rank operators converging to
$A\in \SS_2$ in the $\|\cdot\|_2$ norm, noting that $\SS_2$ has the
approximation property.)

Let then $\lambda_0\in \complex$ be such that $\Det_2(\Id-\lambda_0 A)=0$
with order $m_0\ge1$. For simplicity, we assume that $\lambda_0=1$.
Our assumption is equivalent to the fact that $1$ is an eigenvalue
of algebraic multiplicity $m_0$ for $R_A=\Id-(\Id-A)\exp(A)$.
Let then $\{\varphi_j\, , j=1,\ldots , m_0\}$ be a basis for the
generalised eigenspace of $R_A$ and the eigenvalue $1$.
If $R_A \varphi_j=\varphi_j$ then
$$
-(\Id-A) \exp(A) \varphi_j = - \exp(A) (\Id-A) \varphi_j  = 0 \, ,
$$
so that $\psi_j:= \varphi_j$ is a fixed point of
$A$.

Now, if $\varphi_j$ and $\ell \ge 2$ are such that
$(\Id-R_A)^\ell \varphi_j =0$ but $\varphi'_j=(\Id -R_A)^{\ell-1}
\varphi_j\ne 0$ then
(since $\Id-A$ commutes with $\exp (A)$)
$$
0=\exp (A) (\Id-A) \varphi_j' =
\exp(\ell A) (\Id-A)^\ell  \varphi_j \, ,
$$
while $\exp((\ell-1) A)(\Id-A)^{\ell-1} \varphi_j \ne 0$,
and thus, using commutativity again,
$$ \exp(\ell A)(\Id-A)^{\ell-1}\varphi_j \ne 0\, . $$
Taking $\psi_j:=\varphi_j$, we complete our identification of
the generalised basis of $R_A$  and that
of $A$ for $1$. \qed
\enddemo

\proclaim{Corollary of Exercise 1}
Let $z \mapsto A(z)$ be an analytic map at $z_0 \in \complex$
with $A(z)\in  \SS_2$. Then the map
$z \mapsto \Det_2(\Id+A(z))$  is analytic at $z_0$.
\endproclaim

We also mention for the record:

\proclaim{Corollary of the Lidskii Theorem}
For $A\in \SS_2$, writing $\lambda_j(A)$ for the eigenvalues
of $A$  repeated with multiplicity, we have
$$
\Det_2 (\Id -A) = \prod_j (1-\lambda_j(A))\exp(\lambda_j(A)) \, .
$$
\endproclaim



\enddocument